\newif\ifsmf
\IfFileExists{smfart.cls}{\smftrue}{\smffalse}

\ifsmf
  \documentclass[12pt]{smfart} 
\else
  \documentclass[12pt,a4paper]{amsart}
\fi
\calclayout

\usepackage[T1]{fontenc}
\usepackage[francais]{babel}
\usepackage{mathrsfs}
\usepackage{amssymb}

\theoremstyle{plain}
  \newtheorem{theo}[subsection]{Th\'eor\`eme}
  \newtheorem{prop}[subsection]{Proposition}
  \newtheorem{lemm}[subsection]{Lemme}
  \newtheorem{coro}[subsection]{Corollaire}\theoremstyle{remark}
  \newenvironment{enonce}[1]{\bgroup
     
     \newtheorem{Enonce}[subsection]{#1}\begin{Enonce}}%
     {\end{Enonce}\egroup}
  \newenvironment{enonce*}[1]{\bgroup
     
     \newtheorem*{Enonce}[subsection]{#1}\begin{Enonce}}%
     {\end{Enonce}\egroup}
\theoremstyle{definition}
  \newtheorem{defi}[subsection]{Definition}
  \newtheorem{rema}[subsection]{Remarque}
  \numberwithin{equation}{subsection}
  \def\theequation{\arabic{section}.\arabic{subsection}.\arabic{equation}}
\title[Points de petite hauteur sur les vari\'et\'es semi-ab\'eliennes]
  {Points de petite hauteur \\ sur les vari\'et\'es semi-ab\'eliennes}
\alttitle{Small points on semi-abelian varieties}

\author{Antoine Chambert-Loir}
\address{Institut de math\'ematiques de Jussieu \\
Universit\'e Pierre et Marie Curie (Paris VI) \\
4 place Jussieu, 75252 Paris Cedex 05, France}
\email{chambert@math.jussieu.fr}
\date{Soumis sur l'archive \emph{xxx} le 8 juillet 1999}

\def\gm{\mathbf{G}_{\mathrm{m}}}
\def\norm#1{\left\|{#1}\right\|}
\def\CH{\mathrm{CH}}
\def\hCH^#1{\widehat{\mathrm{CH}}\vphantom{\operatorname{CH}}^{#1}}

\def\hdeg{\widehat{\operatorname{deg}}\,}

\def\Pic{\mathop{{\mathrm{Pic}}}\nolimits}
\def\Arg{\mathop{{\mathrm{Arg}}}\nolimits}

\def\Card{\mathop{{\mathrm{Card}}}\nolimits}
\def\Ext{\operatorname{Ext}}

\def\Cube{\operatorname{Cube}}
\def\Carre{\operatorname{Carr\acute e}}
\def\abs#1{\left|{#1}\right|}
\let\bar\overline
\def\hc{\vphantom{c}{\widehat c}\vphantom{c}}
\def\sozat{\hskip0.3em plus0.2em minus0.1em;\hskip0.5em plus0.3em minus0.1 em}

\let\hat\widehat
\let\tilde\widetilde

\def\C{{\mathbf C}}

\def\N{{\mathbf N}}
\def\P{{\mathbf P}}
\def\Q{{\mathbf Q}}
\def\R{{\mathbf R}}
\def\V{{\mathbf V}}
\def\Z{{\mathbf Z}}

\def\Gal{\operatorname{Gal}}
\def\Spec{\operatorname{Spec}}
\let\union\cup \let\inter\cap

\let\hra\hookrightarrow
\let\ra\rightarrow

\let\eps\varepsilon
\let\phi\varphi\let\mathcal\mathscr
\let\vide\varnothing
\let\leq\leqslant
\let\geq\geqslant

\begin{document}

\begin{abstract}
Dans ce texte, nous \'etendons au cas 
des vari\'et\'es semi-ab\'eliennes les \'enonc\'es
de diverses variantes de la conjecture de Bogomolov
concernant la non-densit\'e
des points de petite hauteur.
Inspir\'es par les travaux r\'ecents 
de Ullmo, Zhang et Bilu, nous d\'emontrons ensuite ces conjectures
lorsque la vari\'et\'e semi-ab\'elienne est isog\`ene \`a un produit
d'une vari\'et\'e ab\'elienne et d'un tore.
La d\'emonstration utilise la g\'eom\'etrie d'Arakelov et des
arguments d'\'equir\'epartition.
\end{abstract}

\ifsmf
\begin{altabstract}
We extend to the case of semi-abelian varieties
the statements of various variants of the conjecture
\emph{alla} Bogomolov about the non-density
of small points of small height in abelian varieties.
Inspired by recent work of Ullmo, Zhang and Bilu,
we then prove these conjectures when the semi-abelian variety
is almost split.
The proof uses Arakelov geometry and equirepartition arguments.
\end{altabstract}
\fi

\maketitle

\section*{Introduction}

Cet article expose quelques propri\'et\'es des hauteurs sur
des vari\'et\'es semi-ab\'eliennes vues sous l'angle de la th\'eorie
d'Arakelov.

Soient $k$ un corps de nombres et $G$ une vari\'et\'e semi-ab\'elienne
(c'est-\`a-dire une extension $0\ra T\ra G\ra A\ra 0$ d'une vari\'et\'e
ab\'elienne~$A$ par un tore~$T$) d\'efinie sur $k$. En tant que
groupe alg\'ebrique, il est quasi-projectif ; on dispose donc
de hauteurs de Weil.

\`A tout triplet $(X,D,\mathcal L)$
constitu\'e d'une compactification \'equivariante lisse $T\hra X$ du tore,
d'un diviseur ample $T$-invariant sur $X$ 
et d'un faisceau inversible, ample et sym\'etrique sur $A$, nous
attachons dans cet article un repr\'esentant canonique de la hauteur de Weil.
Cette {\em hauteur normalis\'ee\/} se comporte bien sous l'action
des morphismes de multiplication par un entier positif ; quand $T=1$,
elle co\"{\i}ncide avec la hauteur de N\'eron--Tate, tandis qu'elle
redonne la {\og hauteur na\"{\i}ve\fg} (de Weil--Mahler) usuelle quand $A=0$
et $X$ est l'espace projectif.
En utilisant la th\'eorie des m\'etriques ad\'eliques de S.~Zhang, nous
d\'efinissons aussi la hauteur canonique d'une sous-vari\'et\'e 
de $G\otimes\bar k$.

Ceci fait, on est en droit de tenter de g\'en\'eraliser aux vari\'et\'es
semi-ab\'eliennes les r\'esultats de Ullmo~\cite{ullmo98} et
Zhang~\cite{zhang98} concernant
la conjecture de Bogomolov.

Le r\'esultat principal de cet article est le suivant
(voir le paragraphe~\ref{sect:conj},
conjectures~\ref{conj:bogo}, \ref{conj:pph}
et~\ref{conj:equirep} et le th\'eor\`eme~\ref{theo:equirep})
o\`u, conform\'ement \`a l'usage,
nous appellerons \emph{sous-vari\'et\'e de torsion}
d'une vari\'et\'e semi-ab\'elienne $G$
d\'efinie sur un corps alg\'ebriquement clos
toute sous-vari\'et\'e qui est la translat\'ee par un point de torsion d'un
sous-groupe alg\'ebrique connexe de~$G$.

\theoremstyle{plain}\newtheorem*{theo*}{Th\'eor\`eme}
\begin{theo*}
Soient $k$ un corps de nombres et
$G$ une vari\'et\'e semi-ab\'elienne d\'efinie sur~$k$
\emph{isotriviale}, c'est-\`a-dire isog\`ene au produit
d'une vari\'et\'e ab\'elienne et d'un tore.
On a les trois \'enonc\'es suivants.

{\normalfont\bfseries\medskip\noindent
a) \'Equir\'epartition des petits points :}
soit $(x_n)_{n\in \N}$ une suite de points de $G(\bar k)$ dont
la hauteur normalis\'ee tend vers $0$ et telle qu'aucune sous-suite
ne soit contenue dans un sous-groupe alg\'ebrique strict.
Alors,
la suite $(\mu_n)_{n\in\N}$ de mesures sur $G(\C)$ d\'efinie par
\[  \mu_n
= \frac{1}{[k(x_n):k]} \sum_{\sigma:k(x_n)\hra\C} \delta_{\sigma(x_n)} \]
(si $x\in G(\C)$, $\delta_x$
d\'esigne comme d'habitude la mesure de Dirac au point~$x$)
converge vaguement vers la mesure de Haar normalis\'ee support\'ee par
le sous-groupe compact maximal de $G(\C)$.

{\normalfont\bfseries\medskip\noindent
b) Sous-vari\'et\'es de hauteur nulle :}
soit $X$ une sous-vari\'et\'e ferm\'ee irr\'eductible de $G\otimes\bar k$.
La hauteur normalis\'ee de $X$ est positive ou nulle ;
elle est nulle si et seulement si
$X$ est une sous-vari\'et\'e de torsion de~$G$.

{\normalfont\bfseries\medskip\noindent
c) Non-densit\'e des petits points :}
soit $X$ une sous-vari\'et\'e ferm\'ee de $G$.
\begin{enumerate}
\item
l'ensemble des sous-vari\'et\'es de torsion de $G$ contenues dans $X$
n'a qu'un nombre fini d'\'el\'ements maximaux ;
\item
notons $X^*$ le compl\'ementaire dans $X$ de la r\'eunion des sous-vari\'et\'es
de torsion contenues dans $X$, 
il existe alors $\eps>0$ tel que l'ensemble
des points de $X^*(\bar k)$ dont la hauteur normalis\'ee est inf\'erieure
\`a $\eps$ est vide.
\end{enumerate}
\end{theo*}

Notons que cette hauteur normalis\'ee n'intervient dans le
premier alin\'ea du th\'eor\`eme ci-dessus
que sous la forme {\og les hauteurs d'une suite de points
tendent vers~$0$\fg}, et on peut v\'erifier que cette notion
ne d\'epend que de~$G$.
De m\^eme pour l'existence d'une minoration de la hauteur
normalis\'ee des points de $X^*$ dans le troisi\`eme alin\'ea.

Le th\'eor\`eme, en tout cas son deuxi\`eme alin\'ea, est \emph{faux}
si l'on omet l'hypoth\`ese que $G$ est isotriviale,
cf.\ le th\'eor\`eme~\ref{theo:faux}.

L'\'equivalence des deuxi\`eme et troisi\`eme alin\'eas du th\'eor\`eme
r\'esulte plus ou moins formellement
des travaux de S.~Zhang~\cite{zhang95,zhang95b}
sur l'amplitude arithm\'etique.
Dans le cas o\`u $T=1$,
L.~Szpiro, E.~Ullmo et S.~Zhang ont prouv\'e dans~\cite{szpiro-u-z97}
que les deux premiers alin\'eas du th\'eor\`eme
(\'equir\'epartition {\em vs.} conjecture de Bogomolov)
sont aussi \'equivalents.

\`A ce propos,
notons que le point 1) de l'alin\'ea c) du th\'eor\`eme n'est pas nouveau :
M.~Hindry l'a en effet d\'emontr\'e (voir~\cite{hindry88})
pour tout groupe alg\'ebrique commutatif.
Cependant, nous n'avons pas besoin d'utiliser son th\'eor\`eme :
sous l'hypoth\`ese que la vari\'et\'e semi-ab\'elienne est isotriviale,
c'est un corollaire de l'approche arakelovienne de ces questions.
Le cas $T=1$ \'etait d\^u \`a M.~Raynaud~\cite{raynaud83c}
tandis que le cas $A=0$
a \'et\'e prouv\'e par M.~Laurent~\cite{laurent84}. 

Lorsque $T=1$, ce th\'eor\`eme a \'et\'e d\'emontr\'e r\'ecemment par
E.~Ullmo~\cite{ullmo98}
pour une courbe dans sa jacobienne, puis presqu'imm\'ediatement g\'en\'eralis\'e
par S.~Zhang~\cite{zhang98} ;  peu apr\`es,
S.~David et P.~Philippon ont d\'emontr\'e dans~\cite{david-p98}
l'existence d'une minoration de la hauteur normalis\'ee des sous-vari\'et\'es
qui ne sont pas des vari\'et\'es de torsion.

Lorsqu'en revanche $A=0$,
le premier alin\'ea du th\'eor\`eme est un th\'eor\`eme
de Y.~Bilu~\cite{bilu97},
alors que le deuxi\`eme alin\'ea a \'et\'e d\'emontr\'e dans
le cas d'une hypersurface par W.~Lawton dans~\cite{lawton77},
puis dans le cas g\'en\'eral par S.~Zhang~\cite{zhang95}
en lien avec le troisi\`eme alin\'ea.

Dans le cas o\`u $G=A\times T$ est un produit, l'\'equir\'epartition
d'une suite de petits points $(x_n=(z_n,t_n))$
est ainsi un th\'eor\`eme une fois projet\'e sur chacun des facteurs $A$
et~$T$.
Il convient d'insister sur le fait que cela n'implique
pas notre th\'eor\`eme qui peut \^etre vu comme une d\'emonstration
de \emph{l'ind\'ependance} des deux variables al\'eatoires 
$(z_n)$ et $(t_n)$.

Nous d\'emontrons ce th\'eor\`eme
en suivant la voie inaugur\'ee par E.~Ullmo, c'est-\`a-dire
\`a l'aide d'arguments d'\'equir\'epartition.
Cependant, les m\'etriques hermitiennes naturelles dans le contexte
des vari\'et\'es semi-ab\'eliennes sortent du cadre classique de la
g\'eom\'etrie d'Arakelov :
\begin{itemize}
\item elles ne sont pas lisses ; cela oblige d'une part \`a d\'efinir les
hauteurs par un proc\'ed\'e limite en recourant \`a la th\'eorie
des m\'etriques ad\'eliques de Zhang, et d'autre part \`a d\'efinir 
divers objets analytiques par approximation via un th\'eor\`eme
d\^u \`a Bedford--Taylor et Demailly qui permet d'effectuer
des produits de courants sous des hypoth\`eses de positivit\'e.
L'utilit\'e de cette th\'eorie en g\'eom\'etrie d'Arakelov est apparue
pour la premi\`ere fois dans les travaux de V.~Maillot~\cite{maillot97}
concernant la g\'eom\'etrie d'Arakelov des vari\'et\'es toriques.

\item les {\og courants de courbure\fg}
qui interviennent sont concentr\'es sur le sous-groupe compact maximal,
ce qui engendre des complications techniques. Pour les contourner,
nous traitons tout d'abord le cas du rang torique~$1$ (voir
le lemme~\ref{lemm:Phi} 
et l'appendice quand $G=\gm$).
Nous utilisons ensuite  un argument de Bilu pour passer au cas g\'en\'eral.
\end{itemize}

Comme il nous a sembl\'e que le cas du groupe multiplicatif plong\'e
naturellement dans $\P^1$ \'eclairait la situation, sans pour autant
contenir les complications li\'ees aux vari\'et\'es ab\'eliennes, nous avons
inclus ce cas dans un appendice ; nous conseillons bien s\^ur au lecteur
de le lire en premier.

Signalons enfin que ce th\'eor\`eme a \'et\'e utilis\'e r\'ecemment par
B.~Poonen~\cite{poonen99}
pour prouver un th\'eor\`eme sur les vari\'et\'es semi-ab\'eliennes
isotriviales qui combine les \'enonc\'es
de type Mordell-Lang (prouv\'e par M.~McQuillan, \cite{mcquillan95})
et Bogomolov.
Ce th\'eor\`eme serait d'ailleurs d\'emontr\'e pour les
vari\'et\'es semi-ab\'eliennes quelconque d\`es que l'on aurait
prouv\'e les alin\'eas a) et c) sans l'hypoth\`ese d'isotrivialit\'e.

\bigskip

Le plan de cet article est le suivant :
les \S\,\ref{sect:prelim} et~\ref{sect:canab}
contiennent des rappels pr\'eliminaires de g\'eom\'etrie
d'Arakelov (particuli\`erement en ce qui concerne la th\'eorie
des m\'etriques ad\'eliques de S.~Zhang) et
des compl\'ements concernant les m\'etrisations des fibr\'es
inversibles sur les vari\'et\'es ab\'eliennes.

La construction des hauteurs canoniques est faite au \S\,\ref{sect:constr}
et dans le \S\,\ref{sect:calcul} qui lui fait suite, nous
calculons la hauteur d'une vari\'et\'e semi-ab\'elienne compactifi\'ee
\`a l'aide du plongement $\gm^t\hra \P^t$.

Nous \'enon\c{c}ons au \S\,\ref{sect:conj} les diff\'erentes variantes
de la conjecture de Bogomolov et pr\'ecisons en d\'etail leurs liens.
Puis, dans les \S\,\ref{sect:generique} et~\ref{sect:equirep},
nous d\'emontrons le th\'eor\`eme principal. 

\bigskip

\emph{Remerciements :} 
je voudrais remercier Daniel Bertrand pour m'avoir incit\'e
\`a r\'efl\'echir \`a ces questions, ainsi que
Ahmed Abbes, Sinnou David, Vincent Maillot et Bjorn Poonen
pour leur int\'er\^et dans ce travail.
Je voudrais aussi remercier Klaus K\"unnemann et Michel Raynaud
pour leur aide \`a propos du lemme~\ref{lemm:klaus}.

\tableofcontents

\section{Pr\'eliminaires de g\'eom\'etrie d'Arakelov}
\label{sect:prelim}

Pour la commodit\'e du lecteur, nous commen\c{c}ons par
rappeler bri\`evement dans ce paragraphe 
quelques r\'esultats de g\'eom\'etrie d'Arakelov
dont nous ferons par la suite un usage constant. 

\subsection{G\'eom\'etrie complexe}
Rappelons quelques notations de g\'eom\'etrie complexe.
Soit $X$ une vari\'et\'e analytique complexe.  On note $\partial$ et
$\bar\partial$ les composantes holomorphes et antiholomorphes
de la diff\'erentielle ext\'erieure $d=\partial+\bar\partial$
sur les formes (voire les courants) sur~$X$. On note
$d^c=(i/4\pi)(\bar\partial-\partial)$, si bien que
$ dd^c=-d^c d = (i/2\pi)\partial\bar\partial$.

On appellera \emph{fibr\'e en droites hermitien} sur $X$ un couple
$\overline{\mathcal L}=(\mathcal L,\norm{\cdot})$
constitu\'e d'un fibr\'e en droites holomorphe sur~$X$ et
d'une m\'etrique hermitienne continue sur $\mathcal L$. 

Si $X$ est lisse, on attache \`a un fibr\'e en droites hermitien un courant
de type $(1,1)$ appel\'e {\em courant
de courbure\/} et d\'efini localement par
\[  c_1(\overline{\mathcal L})=- dd^c \log\norm{\mathsf s}^2,  \]
o\`u $\mathsf s$ est une section locale de $\mathcal L$ ne s'annulant pas
(cela ne d\'epend pas du choix de~$\mathsf s$). Si la m\'etrique hermitienne
sur $\mathcal L$ est lisse, son courant de courbure est en fait
une forme diff\'erentielle, dite {\em forme de courbure.}

On dira qu'un courant $T$ de type~$(p,p)$ est positif si pour tout
choix de formes diff\'erentielles de type~$(1,0)$, $\mathcal C^\infty$ 
\`a support compact sur~$X$,
$(\omega_1,\ldots,\omega_{d-p})$, $d$ d\'esignant la dimension de~$X$,
le courant 
\[ T\wedge (i\alpha_1\wedge \bar\alpha_1) \wedge \cdots
    \wedge (i\alpha_{d-p}\wedge \bar\alpha_{d-p}) \]
est une mesure positive sur~$X$. Quand $T=c_1(\overline{\mathcal L})$,
cela \'equivaut \`a dire que la m\'etrique hermitienne sur $\mathcal L$
s'\'ecrit (dans une trivialisation locale $\iota:\mathcal L\simeq\mathcal O_X$)
$\norm{\mathsf s}^2(x)=|\iota(\mathsf s)(x)|^2 \exp(-u(x))$,
o\`u $u$ est une fonction pluri-sous-harmonique
(psh pour abr\'eger) sur~$X$.

\subsection{Cas des espaces singuliers}
Soit $X$ un espace analytique complexe r\'eduit. Suivant~\cite{zhang95},
on dira qu'une fonction $f$ sur $X$ est $\mathscr C^\infty$
si pour tout morphisme analytique $\phi:U\ra X$ d'un ouvert de $\C^n$ vers $X$,
la compos\'ee $\phi^* f $ est une fonction $\mathscr C^\infty$ sur $U$.
On d\'efinit de m\^eme la notion de forme diff\'erentielle
de classe $\mathscr C^\infty$, de fonction psh sur $X$.

Un fibr\'e en droites hermitien $\overline{\mathscr L}$
sur $X$ sera dit $\mathscr C^\infty$ (resp.\ \`a courbure positive)
si dans toute trivialisation locale, la m\'etrique est donn\'ee par
une fonction $\mathscr C^\infty$ (resp.\ psh).
Cela revient \`a dire que pour tout morphisme
analytique $\phi:U\ra X$ d'un ouvert de $\C^n$ vers $X$, le fibr\'e
en droites hermitien $\phi^*\overline{\mathscr L}$ sur $U$ est
muni d'une m\'etrique lisse (resp.\ \`a courbure positive).

On commetra aussi l'abus de langage consistant \`a \'ecrire
que, $\overline{\mathscr L}$ \'etant un fibr\'e en droites hermitien
sur un espace analytique complexe,
le courant $c_1(\overline{\mathscr L})$ est positif.
Cela signifiera que la m\'etrique est donn\'ee par une fonction continue
et psh.

\subsection{M\'etriques ad\'eliques}
Soit $X$ une vari\'et\'e projective sur un corps de nombres~$k$
et $\mathcal L$ un fibr\'e inversible sur~$X$.

La g\'eom\'etrie d'Arakelov consid\`ere alors
la donn\'ee d'un mod\`ele projectif et plat $\tilde X$
sur l'anneau $\mathfrak o_{k'}$
des entiers d'une extension finie $k'$ de~$k$,
d'un fibr\'e inversible $\widetilde{\mathcal L}$
sur $\tilde X$ dont la restriction \`a la fibre g\'en\'erique est 
une puissance de $\mathcal L$,
et d'une m\'etrique hermitienne $\mathcal C^\infty$
sur $X(\C)$, suppos\'ee invariante par la conjugaison
complexe. Si $w$ est une place finie de~$k'$,
on d\'eduit du mod\`ele entier une norme $w$-adique
sur les fibres de $\mathcal L$ en d\'ecr\'etant que les sections enti\`eres
en~$w$ ont une norme $\leq 1$. 
Si pour toute place $v$ de $k$, les m\'etriques $w$-adiques aux
places~$w$ de $k'$ au-dessus de $v$ sont invariantes par les groupes
de Galois locaux,
une telle donn\'ee d\'etermine 
une collection de normes~$v$-adiques sur $\mathcal L$,
appel\'ee {\em m\'etrique ad\'elique alg\'ebrique} sur $\mathcal L$.

Plus g\'en\'eralement, on dit (\cite{zhang95b}, 1.2)
qu'une collection $(\rho_v)$ de normes~$v$-adiques sur
les fibres de $\mathcal L$ est une {\em m\'etrique ad\'elique} s'il
existe une m\'etrique ad\'elique alg\'ebrique $(\tau_v)$ sur $\mathcal L$
telle que :
\begin{itemize}
\item pour toute place~$v$, le rapport $\tau_v/\rho_v$ est une fonction
continue born\'ee sur $X(\bar k_v)$, invariante sous l'action
de~$\Gal(\bar k_v/k_v)$ ;
\item pour presque toute place non archim\'edienne~$v$, $\tau_v=\rho_v$.
\end{itemize}

Une suite de m\'etriques ad\'eliques $(\rho_{v,n})_n$ converge vers
une m\'etrique ad\'elique $(\rho_v)$ si pour presque tout~$v$,
la suite $\rho_{v,n}=\rho_v$, et si pour tout $v$,
la suite $\log(\rho_{v,n}/\rho_v)$ converge
uniform\'ement vers~$0$. 
 
En g\'eom\'etrie d'Arakelov, les m\'etriques ad\'eliques qui interviennent
sont limites de m\'etriques ad\'eliques alg\'ebriques particuli\`eres :
on dit qu'une m\'etrique ad\'elique alg\'ebrique est~{\em ample}
(resp.\ {\em semi-positive\/}) si,
avec les m\^emes notations que ci-dessus, $\widetilde{\mathcal L}$
est ample sur $\tilde X$
(resp.\ a un degr\'e $\geq 0$ sur toute courbe contenue
dans une fibre verticale de $\tilde X$)
et si la forme de courbure de $\mathcal L$
en les places archim\'ediennes est strictement positive (resp.\ positive
ou nulle). 
Une m\'etrique ad\'elique est dite \emph{ample}
(resp.\ \emph{semi-positive})
si elle est limite de m\'etriques ad\'eliques alg\'ebriques amples
(resp.\ semi-positives).

Le produit tensoriel de deux m\'etriques ad\'eliques d\'efinit encore une
m\'etrique ad\'elique ; de m\^eme la m\'etrique ad\'elique duale est
d\'efinie naturellement.
On dit ainsi qu'une m\'etrique ad\'elique $(\rho_v)$ sur $\mathcal L$ est
{\em int\'egrable} s'il existe deux fibr\'es inversibles $\mathcal L_1$
et $\mathcal L_2$ sur $X$ munis de m\'etriques ad\'eliques semi-positives
$(\rho_{1,v})$ et $(\rho_{2,v})$
tels que $\mathcal L=\mathcal L_1\otimes \mathcal L_2^\vee$
et $\rho_v=\rho_{1,v}\otimes\rho_{2,v}^\vee$
({\em cf.}~\cite{zhang95b}, 1.3).
On notera couramment $\overline{\mathcal L}$ pour d\'esigner la donn\'ee
d'un fibr\'e inversible $\mathcal L$ et d'une m\'etrique ad\'elique
sur $\mathcal L$.

\subsection{Hauteurs}

Ces m\'etriques ad\'eliques permettent de d\'efinir
des nombres d'intersection arithm\'etique
et la hauteur des sous-vari\'et\'es. 
En effet, Zhang \'etend dans~\cite{zhang95b}
la th\'eorie de l'intersection arithm\'etique
de Gillet--Soul\'e (\cite{gillet-s90}, 
voir aussi~\cite{soule-a-b-k92} pour une introduction 
et~\cite{bost-g-s94} pour un expos\'e d\'etaill\'e des hauteurs 
dans ce point de vue). Il prouve en effet  
que si $\overline{\mathcal L}_i$ (pour $1\leq i\leq \dim V+1$)
sont des fibr\'es en droites
sur $X$ munis de m\'etriques ad\'eliques semi-positives
et si $V$ est une sous-vari\'et\'e ferm\'ee de $X$,
les nombres d'intersection arithm\'etique 
\[  (\hc_1(\widetilde{\mathcal L}_1) \cdots
      \hc_1(\widetilde{\mathcal L}_{\dim V+1})|V)  \]
(calcul\'es \`a l'aide d'une d\'esingularisation g\'en\'erique de $V$ et
d'un mod\`ele entier $(\widetilde V,(\widetilde {\mathcal L}_i))$,
comme dans~\cite{bost-g-s94}, 2.3, si ce n'est que nous ne nous int\'eressons
qu'au degr\'e arithm\'etique $\hdeg$
de l'\'el\'ement calcul\'e dans~{\em loc.cit.})
convergent lorsque les m\'etriques ad\'eliques des $\widetilde {\mathcal L}_i$
approchent
celles des $\overline{\mathcal L}_i$ (\cite{zhang95b}, Th.~1.4).
La limite, ind\'ependante de la suite
approximante, est not\'ee naturellement
\[  (\hc_1(\overline{\mathcal L}_1) \cdots
        \hc_1(\overline{\mathcal L}_{\dim V+1})|V).  \]
Si $V$ est fix\'e, cela d\'efinit une application lin\'eaire en chacun
des fibr\'es, ce qui permet d'\'etendre la d\'efinition 
\`a tous les fibr\'es
en droites munis de m\'etriques ad\'eliques int\'egrables.

Lorsque ${\mathcal L}$ est un fibr\'e en droites gros
(ce qui veut dire que $c_1(\mathcal L)^{\dim X}>0$)
muni d'une m\'etrique ad\'elique semi-positive (voire int\'egrable), on d\'efinit
la {\em hauteur de $V$} relativement \`a $\overline{\mathcal L}$ par 
la formule
\[  h_{\overline{\mathcal L}}( V) = \frac{1}{(\dim V+1)\deg_{\mathcal L} V}
(\hc_1(\overline{\mathcal L})^{\dim V+1}|V).  \]
Si $x\in X(k)$ est un point rationnel, on a
\[  h_{\overline{\mathcal L}}(x) =
- \sum_{\text{places $v$ de $k$}} [k_v:\Q_v] \log \norm{\mathsf s}_v(x);
\] 
cette expression \'etant ind\'ependante du choix d'une section non nulle
$\mathsf s\in\mathcal L_x$.
Cette derni\`ere formule a un sens
si $\overline{\mathscr L}$ est un fibr\'e en droites muni d'une m\'etrique
ad\'elique quelconque et la fonction $h_{\overline{\mathscr L}}$
est alors un repr\'esentant de la hauteur de Weil pour $\mathscr L$.

Un th\'eor\`eme fondamental de S.~Zhang relie les hauteurs des points et
celles des sous-vari\'et\'es de~$X$ pour un fibr\'e en droites muni
d'une m\'etrique ad\'elique ample. 
Si $\overline{\mathcal L}$ est un fibr\'e en droites
sur~$X$ muni d'une m\'etrique ad\'elique int\'egrable, on note
\[ 
 \mu_{\overline{\mathcal L}}(X)  
    = \inf_{x\in X(\bar k)} h_{\overline{\mathcal L}}(x)
\qquad\text{et}\qquad
 e_{\overline{\mathcal L}}(X) 
    = \sup_{Y\subsetneq X}
            \inf_{x\not\in Y(\bar k)} h_{\overline{\mathcal L}}(x) .
\] 
Ce sont des \'el\'ements de $\R\union\{-\infty\}$ et
$\R\union\{\pm\infty\}$ respectivement.

\begin{theo}[{{\em cf.}\ \cite{zhang95b}, Th.~1.10}]
\label{theo:1.10}
Soit $\overline{\mathcal L}$
un fibr\'e inversible gros muni d'une m\'etrique ad\'elique semi-positive
et tel que $\mu_{\overline{\mathcal L}}(X)\neq -\infty$.
Alors, on a l'in\'egalit\'e
\[  e_{\overline{\mathcal L}}(X) \geq h_{\overline{\mathcal L}}(X) 
\geq \frac{1}{\dim X+1} \left( e_{\overline{\mathcal L}}(X) +
(\dim X) \mu_{\overline{\mathcal L}}(X) \right).  \]
\end{theo}


\begin{proof}
Ce th\'eor\`eme est prouv\'e dans~\cite{zhang95b}
lorsque la m\'etrique ad\'elique est ample. Nous pourrons ainsi l'appliquer
lorsque le fibr\'e en droites est ample sur la fibre g\'en\'erique (ce qui 
implique que $\mu_{\overline{\mathcal L}}>-\infty$) et muni
d'une m\'etrique ad\'elique semi-positive.

Soit $\pi:\tilde X\ra X$ une d\'esingularisation g\'en\'erique de $X$
et $\overline{\mathcal M}$ un fibr\'e inversible muni d'une
m\'etrique ad\'elique ample sur $\tilde X$.
Pour $a\in\N^*$, soit $\overline{\mathscr M_a}=\pi^*\overline{\mathcal
L}^{\otimes a}\otimes \overline{\mathcal M}$ ;
c'est un fibr\'e en droites ample sur $\tilde X$ muni d'une m\'etrique
ad\'elique semi-positive. D'apr\`es le th\'eor\`eme d'approximation 3.6.1
de~\cite{maillot97}, $\overline{\mathscr M_a}$ est un fibr\'e
inversible muni d'une m\'etrique ad\'elique ample. On peut lui appliquer
le th\'eor\`eme 1.10 de~\cite{zhang95b}, d'o\`u l'in\'egalit\'e
\begin{equation}
\tag{$\dag$}\label{zhang_a}
e_{\overline{\mathscr M_a}} (\tilde X) \geq 
   h_{\overline{\mathscr M_a}} (\tilde X) \geq
   \frac1{1+\dim X} ( e_{\overline{\mathscr M_a}} (\tilde X)
        +  \mu_{\overline{\mathscr M_a}} (\tilde X) ) . 
\end{equation}
Or, en utilisant la multilin\'earit\'e des produits
d'intersection, on voit que lorsque $a$ tend vers $+\infty$,
\begin{enumerate} \def\labelenumi{\theenumi\textsuperscript{o})}
\item $\displaystyle a^{-\dim X} c_1(\pi^*\mathscr L^{\otimes
a}\otimes\mathscr M)^{\dim X}$ converge vers
$c_1(\pi^*\mathscr L)^{\dim X}=c_1(\mathscr L)^{\dim X}$ d'apr\`es
la formule de projection ;
\item de m\^eme, $\displaystyle a^{-1-\dim X} (\hc_1(\pi^*\overline{\mathscr
L}^{\otimes a}\otimes\overline{\mathscr M})^{1+\dim X}|\tilde X)$
converge vers $(\hc_1(\overline{\mathscr L})^{1+\dim X}|X)$ ;
\item comme $\mu_{\overline{\mathscr L}}\neq -\infty$, 
$\displaystyle a^{-1} \mu_{\overline{\mathscr M_a}}$
tend vers $\mu_{\overline{\mathscr L}}$ ;
\item et finalement, puisque $\pi$ est birationnel,
$\displaystyle a^{-1} e_{\overline{\mathscr M_a}}$
tend vers $e_{\overline{\mathscr L}}$.
\end{enumerate}
Lorsqu'on divise les 3 membres de l'in\'egalit\'e~\eqref{zhang_a}
par $a$ et qu'on fait tendre $a$ vers~$+\infty$, 
on obtient l'in\'egalit\'e souhait\'ee.
\end{proof}

\subsection{Th\'eorie de Bedford--Taylor}

En un certain sens, la th\'eorie des m\'etriques ad\'eliques est
une compl\'etion de la th\'eorie d'Arakelov usuelle, en tout cas
en ce qui concerne des produits de premi\`eres classes de Chern
arithm\'etiques. Elle fait intervenir des limites d'int\'egrales
de formes diff\'erentielles positives et le fait remarquable
est que la limite de ces formes existe souvent en tant que courant.

\begin{enonce}{Proposition--D\'efinition}\label{defi:bt}
Soient $\mathcal L_j$ (pour $1\leq j\leq d$) des fibr\'es en droites 
sur une vari\'et\'e projective complexe $X$ de pure dimension~$d$,
munis de m\'etriques
hermitienne continues $\omega_j$ que l'on suppose limites uniformes
de m\'etriques hermitiennes $\omega_{j,k}$ lisses \`a
formes de courbure positives.

Soit $\pi:\tilde X\ra X$ une r\'esolution des singularit\'es de $X$
({\em i.e.}~un morphisme birationnel propre tel que $\tilde X$ est lisse).
Les formes diff\'erentielles
$\bigwedge_{j=1}^d c_1(\pi^*\mathcal L_j,\pi^*\omega_{j,k})$
sur $\tilde X$ convergent
vers un courant positif ferm\'e sur $\tilde X$. Son image par $\pi_*$
est une mesure sur $X$ ind\'ependante du choix de $\tilde X$.
On la note $\bigwedge_{j=1}^d c_1(\mathcal L_j,\omega_j)$ ; cette
expression est lin\'eaire sym\'etrique en les $(\mathcal L_j,\omega_j)$.

Si tous les fibr\'es sont \'egaux \`a $(\mathcal L,\omega)$,
on d\'efinit ainsi $c_1(\mathcal L,\omega)^d$
qui est une mesure positive sur $X$.
Si elle donne \`a $X$ une mesure non nulle,
on notera $\nu_{\mathcal L,\omega}$ l'unique mesure de probabilit\'e
proportionnelle \`a $c_1(\mathcal L,\omega)^d$.
\end{enonce}
\begin{proof}
C'est essentiellement une reformulation des r\'esultats de Bedford--Taylor
et Demailly ({\em cf.}~\cite{bedford-t82}, \cite{demailly93}) sur l'op\'erateur
de Monge--Amp\`ere.
Commen\c{c}ons par prouver le r\'esultat
concernant la limite de courants sur $\tilde X$.
Le probl\`eme \'etant local,
quitte \`a choisir une partition de l'unit\'e finie de $\tilde X$,
il n'est pas restrictif de supposer que les fibr\'es $\mathcal L_j$ sont
triviaux.
Les m\'etriques $\omega_j$ et $\omega_{j,k}$ correspondent alors
\`a des fonctions $u_j$ (continue) et $u_{j,k}$ (lisse)
sur $\tilde X$ via la formule 
\[ \omega_j(f)(x)= |f(x)| e^{-u_j(x)}.  \]
De plus, $u_{k,j}$ converge uniform\'ement vers $u_j$ sur $\tilde X$ et
\[  c_1(\pi^*\mathcal L_j,\pi^*\omega_{j,k}) = dd^c u_{j,k}.  \]
D'apr\`es par exemple le corollaire 2.6 de~\cite{demailly93},
le produit $\bigwedge_{j=1}^d (dd^c u_{j,k})$
converge au sens des courants vers un courant
positif ferm\'e not\'e $\bigwedge_{j=1}^d(dd^c u_j)$ sur $\tilde X$.
Comme il est de degr\'e $(d,d)$, c'est une mesure.

Pour d\'emontrer que l'image de cette mesure sur $X$ est ind\'ependante
du choix de la r\'esolution des singularit\'es choisie, remarquons
qu'on peut {\og coiffer\fg} deux r\'esolutions $\tilde X_1\ra X$
et $\tilde X_2\ra X$  par une troisi\`eme $\tilde X_3\ra \tilde X_1\times_X
\times \tilde X_2 \ra X$. 
Il suffit alors de d\'emontrer que la mesure
calcul\'ee sur $\tilde X_3$ pouss\'ee sur $\tilde X_1$ est \'egale \`a la mesure
calcul\'ee sur $\tilde X_1$. Il est clair qu'elles co\"{\i}ncident
l\`a o\`u $\tilde X_3\ra \tilde X_1$ est un isomorphisme. Le compl\'ementaire
dans $\tilde X_3$
est un ferm\'e alg\'ebrique de dimension $<d$ ; c'est un ensemble
pluripolaire qui est donc de mesure nulle pour la mesure
$(dd^c u_3)^d$, voir par exemple~\cite{klimek91}, proposition~4.6.4.
\end{proof}

Les mesures donn\'ees par la proposition pr\'ec\'edente apparaissent
naturellement en g\'eom\'etrie d'Arakelov. C'est dans le travail
de V.~Maillot~\cite{maillot97} sur la g\'eom\'etrie d'Arakelov
des vari\'et\'es toriques qu'elles ont \'et\'e introduites pour la premi\`ere
fois de mani\`ere g\'en\'erale.

Nous utiliserons le calcul suivant :
\begin{prop} \label{prop:calcul}
Soient $X$ une vari\'et\'e projective
d\'efinie sur un corps de nombres $k\subset\C$
et
$\overline{\mathcal L}$ un fibr\'e en droites (big) sur~$X$
muni d'une m\'etrique ad\'elique semi-positive.
Soit $f:X(\C)\ra\R$ une fonction continue. On suppose qu'il
existe un fibr\'e inversible ample m\'etris\'e $\overline{\mathcal M}$
sur $X$ dont la m\'etrique hermitienne sous-jacente 
a une courbure $c_1(\overline{\mathcal M})$ positive
et tel que $c_1(\overline{\mathcal M})+ dd^c f$ est un courant
positif. D\'efinissons le fibr\'e
inversible m\'etris\'e $\overline{\mathcal L}(\eps f)$ comme
$\overline{\mathcal L}$ mais dont la m\'etrique hermitienne \`a la place
infinie $k\hra\C$ est multipli\'ee par $\exp(-\eps f)$.

Le fibr\'e $\overline{\mathcal L}(\eps f)$ est alors int\'egrable et la hauteur
de $X$ relativement \`a ce fibr\'e v\'erifie le d\'eveloppement limit\'e
\[  h_{\overline{\mathcal L}(\eps f)} = h_{\overline{\mathscr L}}(X)
+ \eps \int_{X(\C)} f\nu_{\overline{\mathcal L}} + O(\eps^2).  \]
\end{prop}
\begin{proof}
L'existence de $\overline{\mathcal M}$ montre que le fibr\'e trivial
$\mathcal O_X(\eps f)$ dont la m\'etrique ad\'elique est multipli\'ee
par $\exp(-\eps f)$ est int\'egrable. Ainsi, par produit tensoriel,
$\overline{\mathcal L}(\eps f)$ est int\'egrable.
Comme le produit d'intersection des fibr\'es en droites munis de m\'etriques
ad\'eliques int\'egrables est multilin\'eaire, on a l'\'egalit\'e
\[  (\hc_1(\overline{\mathcal L}(\eps f))^{\dim X+1}|X)
= \sum_{k=0}^{\dim X+1} \eps^k \binom{\dim X+1}{k}
(\hc_1(\overline{\mathcal O_X}(f))^k
\hc_1(\overline{\mathcal L})^{\dim X+1-k}|X).  \]
Compte tenu de la d\'efinition de $h(X)$,
il suffit, pour \'etablir la proposition, de prouver l'\'egalit\'e
\[  (\hc_1(\overline{\mathcal O_X}(f))
\hc_1(\overline{\mathcal L})^{\dim X}|X)
= \int_{X(\C)} f c_1(\overline{\mathcal L})^{\dim X}.  \]
Mais cette derni\`ere \'egalit\'e est vraie lorsque la m\'etrique ad\'elique
est alg\'ebrique (c'est la formule $a(\eta)x=a(\eta\omega(x))$
de~\cite{soule-a-b-k92}, III, 2.3.1) ;
la hauteur est d\'efinie 
par passage \`a la limite en approchant la m\'etrique ad\'elique
par des m\'etriques alg\'ebriques,
de m\^eme que la mesure dans la d\'efinition~\ref{defi:bt}.
\end{proof}

\section{M\'etriques ad\'eliques canoniques sur les vari\'et\'es ab\'eliennes}
\label{sect:canab}

Une partie des r\'esultats de ce paragraphe est bien connue
des sp\'ecialistes. Les pr\'ecisions qu'on peut y trouver nous seront
n\'eanmoins utiles dans la suite de ce texte.

Soit $A$ une vari\'et\'e ab\'elienne sur un corps $k$.
Si $\eps\in A(k)$ est la section neutre,
un fibr\'e en droites rigidifi\'e sur $A$ est la donn\'ee
d'un fibr\'e en droites $\mathscr L$ sur $A$ et d'un 
\'el\'ement non nul de $\mathscr L|_\eps$ (c'est-\`a-dire d'un isomorphisme
$\eps^*\mathscr L\simeq\mathscr O_{\Spec k}$).
Si l'on se donne deux fibr\'es inversibles rigidifi\'es
isomorphes, il existe un unique
isomorphisme, dit rigidifi\'e, qui respecte les rigidifications.

Si $\mathscr L$ est un fibr\'e inversible sur $A$, on note
$\Cube(\mathscr L)$ le fibr\'e inversible sur $A^3$ donn\'e par
\[  p_{123}^*\mathscr L \otimes p_{12}^*\mathscr L^{-1}
\otimes p_{23}^*\mathscr L^{-1}
\otimes p_{31}^*\mathscr L^{-1}
\otimes p_1^*\mathscr L
\otimes p_2^*\mathscr L
\otimes p_3^*\mathscr L  \]
o\`u $p_I$ d\'esigne le morphisme $A^3\ra A$ donn\'e par la somme
des composantes d'indices contenus dans $I$.
Si $\mathscr L$ est rigidifi\'e, $\Cube(\mathscr L)$
est canoniquement rigidifi\'e et
le th\'eor\`eme du cube implique alors que pour tout fibr\'e en droites
rigidifi\'e, il existe un unique isomorphisme $\Cube(\mathscr
L)\simeq\mathscr O_{A^3}$ compatible aux rigidifications ;
cet isomorphisme sera appel\'e isomorphisme du cube.

De m\^eme, $\Carre(\mathscr L)$ est le fibr\'e inversible sur $A^2$
d\'efini avec des notations analogues par 
\[  p_{12}^*\mathscr L
    \otimes p_1^*\mathscr L^{-1} \otimes p_2^*\mathscr L^{-1}. \]
Lorsque $\mathscr L$ est rigidifi\'e, 
$\Carre(\mathscr L)$ aussi
et si $\mathscr L$ est alg\'ebriquement \'equivalent \`a z\'ero,
il existe un unique isomorphisme rigidifi\'e
$\Carre(\mathscr L)\simeq\mathscr O_{A^2}$,
appel\'e isomorphisme du carr\'e.

\begin{lemm}
Soit $A$ une vari\'et\'e ab\'elienne sur un corps $k$
et $\mathscr L$ un fibr\'e en droites rigidifi\'e sur $A$
muni d'un isomorphisme rigidifi\'e 
\[  \phi : [n]^*\mathscr L \simeq \mathscr L^d  \]
pour deux entiers $n$ et $d\geq 2$.

On suppose que 
\begin{itemize}
\item soit $k$ est un corps valu\'e et $\mathscr L$ est muni d'une m\'etrique
$v$-adique,
\item soit $k$ est un corps de nombres et $\mathscr L$ est muni
d'une m\'etrique ad\'elique ;
\end{itemize}
et on suppose que cette m\'etrique est
telle que l'isomorphisme $\phi$ est une isom\'etrie.

L'isomorphisme du cube est alors une isom\'etrie.

Si $\mathscr L$ est alg\'ebriquement \'equivalent \`a z\'ero, on a la m\^eme
assertion avec l'isomorphisme du carr\'e.
\end{lemm}
\begin{proof}
Qu'un morphisme de fibr\'es en droites m\'etris\'es soit une isom\'etrie
se v\'erifie place par place. Il suffit donc de traiter le cas o\`u
$k$ est un corps valu\'e.

Soit $\mathsf s\in\Gamma(A^3,\Cube(\mathscr L))$ l'image de la section
$1$ par l'isomorphisme du cube.
Alors, $\sigma =[n]^*\mathsf s/\mathsf s^d$ est une section globale
rigidifi\'ee
du fibr\'e inversible rigidifi\'e
$[n]^*\Cube(\mathscr L)\otimes\Cube(\mathscr L)^{-d}$.
Comme $\phi$ induit une isom\'etrie de ce fibr\'e inversible avec
le fibr\'e inversible trivial (muni de la m\'etrique triviale $\norm{1}=1$),
la norme de $\sigma$ est constante, \'egale \`a $1$.
Soit $f(x)=\log \norm{\mathsf s}(x)$. C'est une fonction continue
born\'ee sur $A(\bar k)$ telle que pour tout $x\in A(\bar k)$, 
$f(nx)=d f(x)$.
Comme la multiplication par $n$ est surjective sur $A(\bar k)$, on en
d\'eduit que 
\[ \inf f =  d \inf f, \quad  \sup f = d\sup f. \]
Comme $d\geq 2$, $\inf f = \sup f =0$.
Autrement dit, $\mathsf s$ est de norme constante \'egale \`a $1$
et l'isomorphisme du cube est une isom\'etrie.

L'argument pour l'isomorphisme du carr\'e est identique.
\end{proof}

Pour tout fibr\'e en droite rigidifi\'e
sur une vari\'et\'e ab\'elienne qui est sym\'etrique ou antisym\'etrique,
il existe un unique isomorphisme $[n]^*\mathscr L\simeq\mathscr L^d$
avec $d=n^2$ ou $d=n$ suivant les cas.
Le th\'eor\`eme 2.2(b) de~\cite{zhang95b} implique alors l'existence
de m\'etriques rendant ces isomorphismes des isom\'etries.
Le lemme pr\'ec\'edent implique donc que tous les fibr\'es en droite rigidifi\'es
sur une vari\'et\'e ab\'elienne sur un corps valu\'e (resp.\ sur un corps
de nombres) disposent d'une m\'etrique canonique qui rend
les isomorphismes du cube ou du carr\'e une isom\'etrie.

Si $\mathscr L$ est ample et sym\'etrique, il est bien connu que
$h_{\overline{\mathscr L}}$ est la hauteur de N\'eron--Tate
sur $A(\bar k)$ associ\'ee \`a $\mathscr L$, hauteur que l'on
note aussi $\hat h_{\mathscr L}$.

\begin{coro} \label{coro:Q-L}
Soit $A$ une vari\'et\'e ab\'elienne sur un corps $k$ qui est
soit un corps valu\'e soit un corps de nombres.
Soient $\mathscr L$ et $\mathscr Q$ deux fibr\'es en droites rigidifi\'es
sur $A$, $\mathscr L$ \'etant ample et $\mathscr Q$ alg\'ebriquement
\'equivalent \`a z\'ero.
Soit $q\in A(k)$ tel que, notant $t_q:A\ra A$ la translation par
$q$ dans $A$, on ait
\[ \phi_{\mathscr L}(q)=t_q^* \mathscr L\otimes\mathscr L^{-1}=\mathscr Q .\]
Alors, il existe une isom\'etrie canonique compatible aux rigidifications
\[  \overline{\mathscr Q} \simeq t_q^* \overline{\mathscr L}
 \otimes \overline{\mathscr L}^{-1} \otimes \overline{\mathscr L}^{-1}|_q.
\] 
\end{coro}
\begin{proof}
Notons $\mathscr Q'$ le fibr\'e m\'etris\'e du membre de droite. Il
est rigidifi\'e gr\^ace \`a la pr\'esence du troisi\`eme terme si bien
qu'il suffit de prouver qu'il
v\'erifie le th\'eor\`eme du carr\'e avec m\'etrique.
Or,
\begin{align*}
\Carre(\overline{\mathscr Q'}) &=
  p_{12}^* \overline{\mathscr Q'}
 \otimes p_1^*\overline{\mathscr Q'}^{-1}
 \otimes  p_2^*\overline{\mathscr Q'}^{-1} \\
&= p_{12}^* t_q^* \overline{\mathscr L}
    \otimes p_{12}^* \overline{\mathscr L}^{-1}
\otimes p_1^*t_q^* \overline{\mathscr L}^{-1}
    \otimes p_2^* t_q^* \overline{\mathscr L}^{-1} 
\otimes p_1^* \overline{\mathscr L}
    \otimes p_2^* \overline{\mathscr L}
    \otimes \overline{\mathscr L}|_q \\
&= \Cube(\overline{\mathscr L})|_{A^2\times\{q\}}
= \overline{\mathscr O_{A^2}}.
\qquad\qquad\qed
\end{align*}
\let\qed\relax
\end{proof}

Soit $A$ une vari\'et\'e ab\'elienne sur un corps de nombres $k$.
Si $\mathscr L$ est un fibr\'e inversible rigidifi\'e, sym\'etrique
et ample sur $A$,
la m\'etrique ad\'elique canonique sur $\mathscr L$ est elle-m\^eme ample
en vertu du th\'eor\`eme 2.2 de~\cite{zhang95b}.
Lorsque $\mathscr L$ est alg\'ebriquement \'equivalent \`a z\'ero,
le corollaire pr\'ec\'edent implique que la m\'etrique ad\'elique
canonique est int\'egrable. Le fait qu'elle soit m\^eme semi-positive
d\'epend d'une construction alg\'ebro-g\'eom\'etrique plus pr\'ecise
dont Michel Raynaud et Klaus K\"unnemann m'ont donn\'e deux
d\'emonstrations diff\'erentes. Il faut en effet trouver un mod\`ele
semi-positif (d'une puissance) de $\mathscr L$
sur l'anneau des entiers d'une extension de
$k$. Aux places de bonne r\'eduction, cela ne pose aucun probl\`eme,
de m\^eme aux places infinies o\`u la forme de courbure
de la m\'etrique hermitienne canonique est nulle.
Aux places de mauvaise r\'eduction,
il faut montrer qu'une puissance de $\mathscr L$
poss\`ede un mod\`ele num\'eriquement \'equivalent \`a z\'ero.

\begin{lemm}[K\"unnemann] \label{lemm:klaus}
Soit $R$ un anneau de valuation discr\`ete dont le corps des fractions~$k$
est de caract\'eristique~$0$. Soient $A$ une vari\'et\'e ab\'elienne sur~$k$
et $\mathcal Q\in\Pic^0 A$
un fibr\'e inversible alg\'ebriquement \'equivalent \`a~$0$
sur $A$. Alors, il existe une extension finie~$k'$ de~$k$,
un sch\'ema $\tilde A$ propre et plat sur la cl\^oture int\'egrale~$R'$
de $R$ dans~$k'$ et un fibr\'e en droites $\widetilde{\mathcal Q}$
sur $\tilde A$ tels que
\begin{itemize}
\item
$\tilde A\otimes_{R'}k'=A\otimes_k k'$
et $\widetilde{\mathcal Q}\otimes_{R'}k'$ est une puissance $\geq 1$
de~$\mathcal Q$ ;
\item 
pour toute courbe $C$ contenue dans la fibre sp\'eciale de
$\tilde A$, le degr\'e de $\widetilde{\mathcal Q}$
sur $C$ est \'egal \`a~$0$.
\end{itemize}
\end{lemm}
\begin{proof}
En fait, K.~K\"unnemann d\'emontre dans~\cite{kunnemann96},
lemme~8.1, que ceci est vrai pour toute vari\'et\'e projective lisse sur $k$
qui admet un mod\`ele r\'egulier sur $R$
et tout cycle alg\'ebriquement \'equivalent \`a $0$.
Il d\'emontre dans~\cite{kunnemann98} (th\'eor\`emes 3.5, 4.2)
que les vari\'et\'es ab\'eliennes admettent un mod\`ele r\'egulier
quitte \`a effectuer une extension finie de $k$.
\end{proof}

\begin{rema}
Quelques remarques tout de m\^eme sur ce lemme.
Remarquons pour commencer qu'il est 
vrai lorsque la vari\'et\'e ab\'elienne a potentiellement
bonne r\'eduction : sur une extension~$k'$ telle que $A\otimes_k k'$
a bonne r\'eduction, on choisit pour $\tilde A$ le mod\`ele de N\'eron
de $A\otimes_k k'$ ; c'est un sch\'ema ab\'elien et tout fibr\'e en droites
alg\'ebriquement \'equivalent \`a~$0$ sur sa fibre g\'en\'erique
se prolonge en un
unique fibr\'e en droites alg\'ebriquement \'equivalent \`a~$0$ (existence
et propret\'e du sch\'ema ab\'elien dual).

Il est aussi v\'erifi\'e si la vari\'et\'e ab\'elienne est une courbe
elliptique : on peut choisir  pour $\tilde A$ le mod\`ele minimal
r\'egulier et l'\'etude de la matrice d'intersection de la fibre sp\'eciale
montre que tout diviseur de degr\'e~$0$ s'\'etend en un unique diviseur
(\`a coefficients rationnels) qui est de degr\'e~$0$ sur chaque composante
irr\'eductible de la fibre sp\'eciale.

Dans le cas g\'en\'eral, K\"unnemann se ram\`ene au cas des courbes en utilisant
le fait qu'un fibr\'e en droite alg\'ebriquement \'equivalent \`a $0$
provient via une correspondance d'un fibr\'e en
droite de degr\'e~$0$ sur une courbe projective lisse,
quitte \`a faire une extension des scalaires. Pour les courbes,
le m\^eme argument concernant la matrice d'intersection
de la fibr\'e sp\'eciale fournit l'extension voulue.
L'existence d'un mod\`ele r\'egulier o\`u la th\'eorie de l'intersection de Fulton
puisse s'appliquer permet d'\'etendre cette correspondance
et donc d'obtenir l'extension recherch\'ee.
\end{rema}

Appliqu\'ee \`a une courbe elliptique, la formule de Faltings--Hriljac
exprime l'auto-intersection d'un diviseur d'Arakelov de degr\'e $0$
en fonction de la hauteur de N\'eron--Tate du point correspondant.
Sur une vari\'et\'e ab\'elienne, elle admet la g\'en\'eralisation suivante.

\begin{theo} \label{theo:Q2L}
Soit $A$ une vari\'et\'e ab\'elienne sur un corps de nombres $k$,
$\overline{\mathscr L}$ et $\overline{\mathscr Q}$ deux
fibr\'es en droites rigidifi\'es sur $A$ munis de leurs m\'etriques ad\'eliques
canoniques.
On suppose que $\mathscr L$ est ample sym\'etrique et $\mathscr Q$ 
alg\'ebriquement \'equivalent \`a z\'ero. Soit $q\in A(\bar k)$ un point
tel que
$\mathscr Q = \phi_{\mathscr L}(q)$
(${}=t_q^*\mathscr L\otimes\mathscr L^{-1}$).

Alors, pour tout $a\in\{0,\ldots,1+\dim A\}$, $a\neq 2$, on a
\[  (\hc_1(\overline{\mathscr Q})^a \hc_1(\overline{\mathscr L})^{1+\dim
A-a}|A)
= 0,  \]
tandis que pour $a=2$, on a
\[  (\hc_1(\overline{\mathscr Q})^2 \hc_1(\overline{\mathscr L})^{\dim A-1} | A)
= - \frac{2}{\dim A} c_1(\mathscr L)^{\dim A} \hat h_{\mathscr L}(q).  \]
\end{theo}
\begin{proof}
Notons $g=\dim A$ et soit $f$ la fonction polynomiale d\'efinie par
\begin{equation}
\label{eq:f(x)}
f(x) = \left(\hc_1(\overline{\mathscr L})
            +x\hc_1(\overline{\mathscr Q})\right)^{g+1} 
= \sum_{a=0}^{g+1} \binom{g+1}a x^a \hc_1(\overline{\mathscr Q})^a
\hc_1(\overline{\mathscr L})^{1+\dim A-a}.
\end{equation}
On commence par rappeler que 
\[  \hc_1(\overline{\mathscr L})^{1+g} = 0 \ ; \]
en effet, lorsqu'on remplace $\mathscr L$ par $[2]^*\mathscr L$,
ce nombre est multipli\'e par $4^{1+g}$ puisque $[2]^*\overline{\mathscr L}
\simeq \overline{\mathscr L}^4$ et aussi par
$4^g$ car le degr\'e de la multiplication par $2$ est $4^g$.
(Le m\^eme argument permet de prouver la nullit\'e des expressions du th\'eor\`eme
pour $a\neq 2$.)

Comme $\mathscr Q=\phi_{\mathscr L}(q)$, le corollaire~\ref{coro:Q-L}
implique que  pour $x\in\Z$,
\[  x\hc_1(\overline{\mathscr Q}) = \hc_1( t_{xq}^*\overline{\mathscr L})
        - \hc_1(\overline{\mathscr L})
        - \pi^*\hc_1(\overline {\mathscr L}|_{xq})  \]
o\`u $\pi:A\ra\Spec k$ est le morphisme structural canonique.
Ainsi, si $x$ est entier,
\begin{align*}
f(x) &= \left(  \hc_1(t_{xq}^*\overline{\mathscr L})
            -\pi^* \hc_1(\overline{\mathscr L}|_{xq}) \right)^{g+1} \\
&= \sum_{a=0}^{g+1} \binom{g+1}a (-1)^a
          \pi_* \left( \hc_1(t_{xq}^* \overline{\mathscr
L})^{g+1-a}\right)
\hc_1(\overline{\mathscr L}|_{xq})^{a}
\end{align*}
d'apr\`es la formule de projection.
Or, $\hc_1(\overline{\mathscr L})^{g+1-a}$ s'interpr\`ete comme un
\'el\'ement de $\hCH^{g+1-a}(A)$ 
si bien que
$\pi_*\hc_1(\overline{\mathscr L})^{g+1-a}$ est nul si $a\neq 0$
et si $a\neq 1$.
Pour $a=0$, on a
\[
\hdeg( \hc_1(t_{xq}^* \overline{\mathscr L})^{g+1}) 
 = \left( \hc_1(t_{xq}^* \overline{\mathscr L})^{g+1} | A \right)
 = \left( \hc_1(\overline{\mathscr L})^{g+1} | t_{xq\,*} A \right) 
 = \left( \hc_1(\overline{\mathscr L})^{g+1} | A \right) = 0 .
\]
Pour $a=1$, 
$\pi_*\hc_1(\overline{\mathscr L})^{g}$ appartient \`a 
$\hCH^0(\Spec\mathfrak o_k)$ et s'identifie \`a $c_1(\mathscr L)^g$.
Par suite,
\begin{align}
\notag
f(x) & = - \binom{g+1}1 c_1(\mathscr L)^g \hdeg \hc_1(\overline{\mathscr
L}|_{xq}) \\
\notag
& = - (g+1) c_1(\mathscr L)^g \hat h_{\overline{\mathscr L}}(xq) \\
\label{eq:f(x)2}
& = - (g+1) c_1(\mathscr L)^g \hat h_{\overline{\mathscr L}}(q)\cdot x^2. 
\end{align}
En comparant les formules~\eqref{eq:f(x)} et~\eqref{eq:f(x)2},
on obtient le th\'eor\`eme.
\end{proof}

\section{Hauteurs normalis\'ees sur les vari\'et\'es semi-ab\'eliennes}
\label{sect:constr}

\subsection{Compactifications \'equivariantes}

On consid\`ere un sch\'ema semi-ab\'elien
\[  1 \ra T \ra G \ra A \ra 0  \]
sur un sch\'ema $S$ (disons, noeth\'erien, s\'epar\'e ; rapidement,
$S$ sera le spectre de l'anneau des entiers d'un corps de nombres).
Fixons une compactification \'equivariante $T\hra X$ du tore
(ainsi, $X$ est une $S$-vari\'et\'e torique). On supposera
toujours que $X$ est {\em projective lisse\/} sur $S$.
Alors, l'image directe de $G$ par le morphisme $T\ra X$
d\'efinit un $S$-sch\'ema $\bar G$, ainsi qu'un morphisme
$\pi:\bar G \ra A$ propre (c'est une fibration dont les fibres sont,
localement sur $A$ pour la topologie de Zariski, isomorphes \`a $X$) ;
{\em cf.}~par exemple~\cite{bosch-l-r90}, 10.2, Lemma~6. 

De mani\`ere plus \'el\'ementaire, 
dans le cas o\`u le tore $T$ est d\'eploy\'e sur $S$,
un isomorphisme $T\simeq \gm^t$ nous fournit $t$ fibr\'es en droites
sur $A$, alg\'ebriquement \'equivalents \`a $0$, not\'es
$\mathcal L_1$, \dots, $\mathcal L_t$. Alors, si
$T\hra X$ est donn\'e par un \'eventail (r\'egulier) $\mathfrak S$ de $\R^t$, on
peut d\'ecrire $\bar G$ en recollant des sous-sch\'emas de fibr\'es
vectoriels sur $A$ de la forme 
$\V(\bigoplus\limits_{i=1}^t
    (\bigotimes\limits_{j=1}^t \mathcal L_j^{\otimes n_{ij}}))$.

\begin{rema}
Dans ma th\`ese~\cite{chambert-loir99}, j'ai trait\'e le cas o\`u $T=\gm$. 
Lorsque $S$ est le spectre d'un corps de caract\'eristique~$0$,
P.~Vojta fait cette construction dans son article~\cite{vojta99}.
Auparavant, J-P.~Serre~\cite{serre79}
avait consid\'er\'e la compactification $\gm^t\hra (\P^1)^t$.
Voir aussi Knop--Lange~\cite{knop-l85},
toujours lorsque $S$ est le spectre d'un corps.
Dans un travail r\'ecent avec Yuri Tschinkel
~\cite{chambert-loir-t99},
nous introduisons la notion de torseur arithm\'etique qui \'eclaire
particuli\`erement ces constructions de hauteurs dans des situations
fibr\'ees.
\end{rema}

L'action de $T$ sur $X$ se prolonge en une action $m:G\times \bar G\ra\bar G$.
On dispose aussi de $S$-morphismes $\bar G\ra\bar G$,
not\'es $[n]_{\bar G}$ qui prolongent
la multiplication par $n$ : $G\ra G$.

\subsection{Fibr\'es en droites $G$-lin\'earis\'es}


Comme $X$ est lisse, tout diviseur de $X$ invariant par l'action de $T$ 
d\'efinit un unique faisceau inversible $\mathcal L$ sur $\bar G$, 
$G$-lin\'earis\'e au sens que si $p_1$, $p_2$ (resp. $m$) d\'esignent
les deux projections (resp.\ la loi d'addition) $G\times_S \bar G\ra  \bar G$,
on a un isomorphisme
$(p_1,m)^*\mathcal L \simeq (p_1,p_2)^*\mathcal L$.
(Autrement dit, on dispose d'isomorphismes compatibles, pour
$x\in \bar G$ et $g\in G$,
$\mathcal L_{g\cdot x} \simeq \mathcal L_g\otimes \mathcal L_x$.)
En particulier, la restriction de $\mathcal L$ \`a la section
unit\'e de $G$ est munie d'une trivialisation canonique.
De m\^eme, on dispose d'un isomorphisme canonique
$[n]_{\bar G}^*\mathcal L\simeq \mathcal L^{\otimes n}$.

D\'esormais, on se place sur le spectre d'un corps de nombres $k$.

\begin{prop} \label{prop:adelique}
Tout faisceau inversible $G$-lin\'earis\'e $\mathcal L$ sur la
compactification $\bar G$
poss\`ede une unique m\'etrique ad\'elique telle que pour tout entier $n>1$,
l'isomorphisme canonique
$[n]_{\bar G}^*\mathcal L\simeq \mathcal L^{\otimes n}$ soit une isom\'etrie.
\end{prop}
\begin{proof}
Fixons un entier $n>1$.
S.\ Zhang a montr\'e (\cite{zhang95b}, Theorem~2.2) que $\mathcal L$
poss\`ede pour toute place $v$ de $k$ une unique m\'etrique $v$-adique
continue et born\'ee $\norm{\cdot}_v$ telle que l'isomorphisme canonique
$[n]_{\bar G}^* \mathcal L\simeq \mathcal L^{\otimes n}$
soit une isom\'etrie en~$v$.

Il faut montrer que la collection des m\'etriques
$\norm{\cdot}=(\norm{\cdot}_v)_v$ est une m\'etrique ad\'elique.
Il suffit de v\'erifier que la collection de m\'etriques obtenue
provient pour presque toute place $v$ d'un mod\`ele entier 
sur un ouvert $U$ du spectre de l'anneau $\mathfrak o_k$ des entiers de $k$.
Pour cela, il suffit de d\'emontrer les quatre points suivants :
\begin{itemize}
\item la compactification \'equivariante $\bar G$ se prolonge 
en un sch\'ema propre et plat $\overline{G_U}$ sur~$U$ ;
\item le morphisme  $[n]_{\overline G}$ se prolonge en un morphisme
$\phi:\overline{G_U}\ra \overline{G_U}$ ;
\item le faisceau inversible $\mathcal L$ se prolonge 
en un faisceau inversible $\mathcal L_U$ sur $\overline{G_U}$ ;
\item l'isomorphisme $[n]^* \mathcal L\simeq \mathcal L^{\otimes n}$
se prolonge en un isomorphisme $\phi^*\mathcal L_U\simeq\mathcal
L_U^{\otimes n}$.
\end{itemize}
En effet, ce mod\`ele entier d\'efinit pour toute place~$v$ de $k$
qui domine $U$ une m\'etrique $v$-adique canonique ; les propri\'et\'es
ci-dessus impliquent que cette m\'etrique $v$-adique rend
l'isomorphisme $[n]^* \mathcal L\simeq \mathcal L^{\otimes n}$
une isom\'etrie, ainsi qu'il fallait prouver.

Il existe alors un ouvert~$U$ de $\Spec \mathfrak o_k$ tel que :
\begin{itemize}
\item la $k$-vari\'et\'e ab\'elienne $A$ se prolonge en un sch\'ema ab\'elien
$A_U$ sur $U$ :
\item le $k$-tore $T$ se prolonge en un tore $T_U$ sur $U$ ;
\item la vari\'et\'e torique (lisse) $T\hra X$ se prolonge en
une $U$-vari\'et\'e torique lisse $T_U\hra X_U$.
\end{itemize}
Alors, l'extension $G\in\Ext^1_{\Spec k}(A,T)$ se prolonge en une unique
extension $G_U$ de $A_U$ par $T_U$. Le r\'esultat
est bien connu si $T=\gm$ (existence du sch\'ema ab\'elien dual) ;
il s'en d\'eduit imm\'ediatement pour un tore d\'eploy\'e
puis par descente au cas d'un tore quelconque.
Puis, la compactification \'equivariante~$\overline G$
se prolonge en une compactification
\'equivariante~$\overline{G_U}$ du $U$-sch\'ema semi-ab\'elien~$G_U$.

On peut de plus supposer que $\mathcal L$ provient d'un diviseur
effectif $T$-invariant sur $X$ ; alors, ce diviseur s'\'etend
en un unique diviseur $T_U$-invariant sur $X_U$, puis
en un unique diviseur $G_U$-invariant sur $\overline{G_U}$. Ainsi,
$\mathcal L$ se prolonge en un unique faisceau inversible $G_U$-lin\'earis\'e
$\mathcal L_U$ sur $\overline{G_U}$. On dispose ainsi d'un unique
isomorphisme $[n]_{\overline{G_U}}^* \mathcal L_U
\simeq \mathcal L_U^{\otimes n}$ qui prolonge l'isomorphisme
analogue sur $\Spec k$. 

Montrons maintenant que ces m\'etriques ad\'eliques ne d\'ependent pas du choix
de l'entier $n$. Par unicit\'e, la m\'etrique ad\'elique associ\'ee \`a un
entier $n$ est \'egale \`a la m\'etrique ad\'elique associ\'ee
\`a tout multiple $mn$,
laquelle co\"{\i}ncide avec la m\'etrique ad\'elique pour l'entier $m$,
d'o\`u l'assertion d'unicit\'e.
\end{proof}

\begin{rema}\label{rema:exemple}
Dans le cas o\`u $A=0$ et $X=\P^t$,
cette m\'etrique est une des m\'etriques usuelles sur $\mathcal O_{\P}(1)$ :
si $P\in k(x_0,\ldots,x_n)$ est homog\`ene de degr\'e~$1$,
\[  \norm{P}(x_0:\cdots:x_t) 
    = \frac{|P(x_0,\ldots,x_n)|}{\max(|x_0|,\ldots,|x_t|)}.  \]
Plus g\'en\'eralement, toujours quand $A=0$, elle a \'et\'e consid\'er\'ee
de mani\`ere approfondie par Maillot~\cite{maillot97}
dans son \'etude de la g\'eom\'etrie d'Arakelov des vari\'et\'es toriques. 
Quand $A$ est quelconque et $t=1$, donc $X=\P^1$, je 
l'ai utilis\'ee dans~\cite{chambert-loir99}.
\end{rema}

\begin{prop}
Si $\mathcal L$ provient d'un diviseur invariant ample
sur~$X$, cette m\'etrique ad\'elique est de plus semi-positive. 
En particulier,
les m\'etriques ad\'eliques fournies par la proposition~\ref{prop:adelique}
sont int\'egrables.
\end{prop}
\begin{proof}
Quitte \`a faire une premi\`ere extension des scalaires, on peut
supposer que le tore sous-jacent \`a $G$ est d\'eploy\'e.
Soient alors $\mathcal Q_i$ des fibr\'es en droites alg\'ebriquement \'equivalents
\`a~$0$ sur $A$ qui d\'efinissent la vari\'et\'e semi-ab\'elienne~$G$.
D'apr\`es le lemme~\ref{lemm:klaus}, il existe
une extension finie $k'$ de $k$, des mod\`eles~$\tilde A_i$ de $A$
sur $\Spec\mathfrak o_{k'}$, un entier $N\geq 1$
et des mod\`eles $\widetilde{\mathcal Q}_i$ des $\mathcal Q_i^{\otimes N}$
sur $\tilde A_i$ qui sont de degr\'e~$0$ sur toute courbe contenue
dans une fibre verticale.
Quitte \`a remplacer~$A$
par le produit fibr\'e des~$\tilde A_i$ au-dessus de $\Spec \mathfrak o_{k'}$,
nous pouvons supposer que tous les $\tilde A_i$ sont \'egaux \`a $\tilde A$.

Soit $G'$ la vari\'et\'e semi-ab\'elienne image de $G$ par l'image
directe via le morphisme
$T\xrightarrow{[N]}T$. Par le $A$-morphisme naturel
$\overline{G}\ra\overline{G'}$ d'\'el\'evation \`a la puissance $N$ dans
les fibres, $\mathcal L^{\otimes N}$ s'identifie isom\'etriquement
\`a l'image r\'eciproque d'un fibr\'e $\mathcal L'$ sur $\overline {G'}$.
Il suffit ainsi de montrer que le m\'etrique ad\'elique sur $\mathcal L'$
est semi-positive.

Soit $\tilde G$ le $\mathfrak o_{k'}$-mod\`ele de $\bar G'$ obtenu
en recollant des sous-sch\'emas de fibr\'es vectoriels sur $\tilde A$
de la forme $\V(\oplus_i (\otimes_j \widetilde Q_j^{n_{i,j}}))$.
Soit $\tilde{\mathscr L}$ le mod\`ele entier de $\mathscr L'$ induit
par le m\^eme diviseur $T$-invariant $D_X$ sur $X$
que celui qui d\'efinit $\mathscr L$.
D'apr\`es la th\'eorie des vari\'et\'es toriques, il existe alors
des diviseurs $D_s$ $T$-invariants sur $X$ lin\'eairement \'equivalents
\`a un multiple $n_s D_X$ tels que $\bigcap_s D_s=\vide$.
Sur $\tilde G$, on dispose de diviseurs $[D_s]$ induits
par les $D_s$ dont l'intersection est vide.
De plus, il existe  pour tout $s$
un fibr\'e en droites $\mathscr Q_s$ sur $\mathscr A$
de la forme $\bigotimes_j \widetilde Q_j^{m_{j}}$
tel que $\mathscr L'=\mathscr O([D_s])\otimes\pi^*\mathscr Q_s$.

Soit $C$ une courbe contenue dans une fibre verticale de
$\tilde {G}$. Soit $s$ tel que $C\not\subset [D_s]$.
Alors $\deg \mathscr O([D_s])|_C\geq 0$ et puisque les
$\widetilde Q_j$ sont de degr\'e~$0$ sur l'image $\pi(C)$ de
$C$ dans $\tilde A$, on a
\[ \deg \tilde{\mathscr L}|_C \geq 0. \]

De plus, $\mathcal L'$ poss\`ede une m\'etrique
hermitienne dont la forme de courbure est positive ou nulle : 
il suffit de noter que du point de vue diff\'erentiel hermitien, la
fibration $\bar G \ra A$ est localement un produit, un fibr\'e
en droites alg\'ebriquement \'equivalent \`a~$0$
sur une vari\'et\'e ab\'elienne \'etant plat.
On recopie alors la preuve de~\cite[2.3]{zhang95b} en rempla\c{c}ant
le mot {\og ample\fg} par {\og semi-positif\fg}.
\end{proof}

\subsection{Construction des hauteurs normalis\'ees}

Consid\'erons un fibr\'e inversible $G$-lin\'earis\'e
$\mathcal L$ sur $\bar G$, relativement ample sur~$A$,
muni de sa m\'etrique ad\'elique canonique.
Alors, pour tout fibr\'e inversible ample sym\'etrique $\mathcal M$
sur $A$, $\mathcal L\otimes\pi^*\mathcal M$ est ample.
(C'est une cons\'equence du fait que les fibres de $\bar G\ra A$
sont toutes isomorphes, voir~\cite{vojta96}, Lemma 3.1.)
En munissant $\mathcal M$ de la m\'etrique ad\'elique d\'efinissant
la hauteur de N\'eron--Tate sur $A$ comme dans la
paragraphe~\ref{sect:canab},
on obtient une m\'etrique ad\'elique semi-positive canonique sur
$\mathcal L\otimes\pi^*\mathcal M$, et donc
une fonction hauteur sur les sous-vari\'et\'es de $\bar G$. 
Comme $\mathcal L\otimes\pi^*\mathcal M$ est ample, sa restriction
aux points est un repr\'esentant d'une hauteur de Weil sur $\bar G(\bar k)$.

\begin{defi}\label{defi:h-norm}
Cette fonction, attach\'ee \`a
$(T\hra X,\mathcal L,\mathcal M)$,
sera appel\'ee {\em hauteur normalis\'ee}. On la note $h_{\overline{\mathscr
L}\otimes\pi^*\overline{\mathscr M}}$ et $\hat h$ lorsqu'il n'y
a pas d'ambigu\"{\i}t\'e possible.
\end{defi}

Lorsque la vari\'et\'e torique $T\hra X$ choisie pour compactifier
$G$ est une vari\'et\'e de Fano, on peut choisir pour $\mathcal L$ le 
fibr\'e inversible associ\'e au compl\'ementaire de
$T$ dans $X$ qui est un diviseur $T$-invariant relativement ample
(diviseur anticanonique relatif). Dans ce cas,
la hauteur normalis\'ee
ne d\'epend plus que d'une compactification \'equivariante du tore
et d'un faisceau inversible ample sym\'etrique sur la vari\'et\'e ab\'elienne.
Les deux exemples classiquement utilis\'es ($\gm^g\subset(\P^1)^g$
et $\gm^g \subset \P^g$) sont d'ailleurs dans ce cas.

Il r\'esulte de la d\'efinition des m\'etriques ad\'eliques canoniques
que la hauteur normalis\'ee d'un point $x\in \bar G(\bar k)$ relativement
\`a $\overline{\mathcal L}\otimes\pi^*\overline{\mathcal M}$ v\'erifie
pour tout point $x\in G(\bar k)$ et tout entier~$n\geq 1$ la relation
\[  h_{\overline{\mathcal L}\otimes\pi^*\overline{\mathcal M}}([n]_G x)
= h_{\overline{\mathcal L}}([n]_G x)+h_{\overline{\mathcal M}}(\pi([n]_G x))
= n h_{\overline{\mathcal L}}(x) + n^2 h_{\overline{\mathcal M}}(\pi(x)).
\] 
On a alors le lemme suivant :
\begin{lemm} \label{lemm:positivite}
Pour tout $x\in G(\bar k)$, on a $h_{\overline{\mathscr L}}(x)\geq 0$.
Ainsi, pour tout $x\in G(\bar k)$, $\hat h(x)\geq 0$ ;
on a de plus \'egalit\'e si et seulement si $x$ est de torsion.
\end{lemm}
On notera que ce lemme
implique que pour toute sous-vari\'et\'e $X$ de $\bar G$ rencontrant
$G$, on a l'\'egalit\'e
$e_{\overline{\mathcal L}\otimes\pi^*\overline{\mathscr M}}(X)\geq 0$.
\begin{proof}
Comme $\mathcal L$ est repr\'esent\'e par un diviseur contenu dans
$\bar G\setminus G$, toute hauteur de Weil attach\'ee \`a $\mathcal L$
est minor\'ee sur $G(\bar k)$. La propri\'et\'e de normalisation implique
alors que pour tout $x\in G(\bar k)$, $h_{\overline{\mathcal L}}(x)\geq 0$.

Comme $h_{\overline{\mathscr M}}$ est la hauteur de N\'eron--Tate
sur $A$, elle est positive, ce qui implique $\hat h(x)\geq 0$.
Si de plus $x\in  G(\bar k)$ est de hauteur normalis\'ee
nulle, on voit que $h_{\overline{\mathscr M}}(\pi(x))=h_{\overline{\mathscr
L}}(x)=0$.
Par suite, tous les it\'er\'es $x$, $[2]x$, etc.\ sont de hauteur nulle.
Comme $\mathscr L\otimes\pi^*\mathscr M$ est ample
et que $\hat h$ est un repr\'esentant de la hauteur de Weil
pour ce fibr\'e, le th\'eor\`eme de Northcott
implique qu'ils forment un ensemble fini.
Cela implique que $x$ est de torsion. La r\'eciproque est claire.
\end{proof}

\begin{rema}
Lorsque $t=1$ et $X=\P^1$,
D.\ Bertrand donne dans~\cite{bertrand96} une d\'ecomposition de la hauteur
relative d'un point comme une somme de hauteurs locales, toutes positives
ou nulles. En g\'en\'eral, si $\mathscr L$ est repr\'esent\'e par un diviseur
$G$-invariant effectif dont le support ne rencontre pas $G$, on peut
v\'erifier que sa section canonique est de norme
inf\'erieure ou \'egale \`a $1$ en toute place.
\end{rema}

La proposition suivante r\'esume les propri\'et\'es des hauteurs normalis\'ees
vis \`a vis des morphismes de multiplication par un entier positif.
\begin{prop} \label{prop:normalisation}
Soient $V$ une sous-vari\'et\'e de~$\bar G$ et $n$ un entier $\geq 2$.
\begin{enumerate}
\item
Si $[n]_*V$ d\'esigne le {\em cycle\/} image directe de $V$ par
le morphisme de multiplication par~$n$, on a pour tous entiers
$a$ et $b$ positifs tels que $a+b=\dim V+1$ l'\'egalit\'e
\[  (\hc_1(\overline{\mathcal L})^a\hc_1(\pi^*\overline{\mathcal M})^b|[n]_*V)
= n^{a+2b} (\hc_1(\overline{\mathcal L})^a\hc_1(\pi^*\overline{\mathcal
M})^b|V).  \]
\item
Si $\xi$ est un point de torsion de~$G$ et $V$ une sous-vari\'et\'e
de $\bar G$, alors $V$ et $\xi+V$ ont m\^eme hauteur normalis\'ee.
\end{enumerate}
\end{prop}
\begin{proof}
\begin{enumerate}
\item C'est la formule de projection, {\em cf.}~\cite[2.3.1, (iv)]{bost-g-s94}
dans le cas de m\'etriques ad\'eliques alg\'ebriques ; le cas g\'en\'eral s'en
d\'eduit par approximation.
\item
Soit $n$ un entier $\geq 2$ tel que $[n]_G\xi=0$. 
Comme $[n]_*V=[n]_*(\xi+V)$, la formule pr\'ec\'edente implique
l'\'egalit\'e des degr\'es arithm\'etiques.
Comme $V$ et $\xi+V$ ont m\^eme degr\'e g\'eom\'etrique, l'\'egalit\'e des hauteurs
en r\'esulte.  \qed
\end{enumerate}
\let\qed\relax\end{proof}

\section{Un exemple}
\label{sect:calcul}

Dans ce paragraphe, on on calcule explicitement
quelques uns des invariants introduits plus haut
pour une vari\'et\'e semi-ab\'elienne compactifi\'ee
\`a l'aide du plongement $\gm^t\subset\P^t$.

On se fixe ainsi un corps de nombres $k$, une vari\'et\'e
ab\'elienne $A$ de dimension~$g$ d\'efinie sur $k$, un fibr\'e inversible ample
et sym\'etrique $\mathscr M$ sur $A$ et $t$ \'el\'ements $q_1,\ldots,q_t$
de $A(k)$. Pour $i\in\{1,\ldots,t\}$, on pose
$\mathscr Q_i=\phi_{\mathscr M}(q_i)=t_{q_i}^*\mathscr M\otimes\mathscr
M^{-1}$.
Tous ces fibr\'es en droites sont munis de leur m\'etrique ad\'elique canonique.

Soit $G$ la vari\'et\'e semi-ab\'elienne
\[ G := \prod_{i=1}^t
      \big( \V(\mathscr Q_i^\vee)\setminus\{0\} \big). \]
(Le produit est relativement \`a $A$.
Un point de $G$ est la donn\'ee d'un point $a\in A$ et pour tout $i$,
d'un \'element non nul de $\mathscr Q_i|_a$.)

En explicitant la construction de $\bar G$ faite au \S\,\ref{sect:constr},
on peut identifier $\bar G$ au fibr\'e projectif
\[ \bar G = \P \big(\mathscr O_A \oplus \mathscr Q_1^\vee \oplus\cdots\oplus
\mathscr Q_t^\vee \big) \]
dont un point est la donn\'ee d'un point $a\in A$ et d'\'el\'ements
non tous nuls de $\mathscr O_A|_a,\mathscr Q_1|_a,\ldots,\mathscr Q_t|_a$.
Notons $\pi:\bar G\ra A$ la projection canonique.

Les $t+1$ diviseurs $\gm^t$-invariants de $\P^t$ 
d\'efinis par la nullit\'e des coordonn\'ees
fournissent des diviseurs $D_0,\ldots,D_t$ sur $\bar G$.
Les fibr\'es inversibles correspondants seront not\'es
$\mathscr L_0,\ldots,\mathscr L_t$.
De plus, on a pour tout $i$ l'\'egalit\'e $D_i-D_0\sim c_1(\pi^*\mathscr Q_i)$
dans $\CH^1(\bar G)$
(cf.~\cite[Chap. V, Prop. 2.6]{hartshorne77}).
Si les $\mathscr L_i$ sont munis de leur m\'etrique ad\'elique canonique,
il est ais\'e de v\'erifier que l'\'egalit\'e pr\'ec\'edente est vraie avec m\'etriques :
pour tout $i\in\{1,\ldots,t\}$, 
\[ \hc_1(\overline{\mathscr L_i}) = \hc_1(\overline{\mathscr L_0})
+ \pi^* \hc_1(\overline{\mathscr Q_i}). \]
On notera $\hat{D_i}=\hc_1(\overline{\mathscr L_i})$.
Posons aussi
\[  \overline{\mathscr Q} = \bigotimes_{i=1}^t \overline{\mathscr Q_i},
\quad q=\sum_{i=1}^t q_i
\quad\text{et}\quad
    \overline{\mathscr L}= \bigotimes_i \overline{\mathscr L_i}
   = \overline{\mathscr L_0}^{t+1} \otimes
   \pi^*\overline{\mathscr Q}. \]

\begin{prop}
Pour toute famille de fibr\'es en droites $\overline{\mathscr
F_j}$ ($1\leq j\leq g$) munis de m\'etriques ad\'eliques int\'egrables sur $\bar G$,
on a
\[ \widehat{D_0} \cdot \widehat{D_1} \cdots \widehat{D_t} \cdot
   \prod_{i=1}^t \hc_1(\overline{\mathscr F_i}) = 0 .\]
\end{prop}
\begin{proof}
Choisissons des mod\`eles $\tilde A$, $\tilde{\mathscr Q_i}$,
et $\tilde G=\P (\mathscr O_{\tilde A}
\oplus \tilde{\mathscr Q_1}\vphantom{\tilde{\mathscr Q_1}}^\vee\oplus\cdots)$
sur l'anneau $\mathfrak o$ des entiers de $k$.
Aux places archim\'ediennes, munissons les $\mathscr Q_i$ de leur m\'etriques plate 
canonique.

Le cycle arithm\'etique $\tilde D_i$ d\'efini par l'annulation de la
$i$\textsuperscript{e} coordonn\'ee prolonge $D_i$
et on peut le munir d'un courant de Green canonique
en munissant $\mathscr O_A\oplus \mathscr Q_1\oplus\cdots\oplus\mathscr
Q_t$ de la m\'etrique hermitienne somme directe orthogonale.

Comme les diviseurs $\tilde D_i$ s'intersectent proprement,
$\tilde D_0\cdots\tilde D_t$ est repr\'esent\'e
par un couple $(Z,g_Z)\in \widehat{\mathrm Z}\vphantom{Z}^{t+1}(\tilde G)$
o\`u $Z\in\mathrm Z^{t+1}(\tilde G)$ est un cycle support\'e par l'intersection
des $\tilde D_i$ et $g_Z$ un courant de Green pour $Z$.
Or, l'intersection des $\tilde D_i$ est vide, si bien que $Z=0$.
L'\'equation $dd^c g_Z+\delta_Z=\omega_Z$, le fait que $\omega_Z$
soit $\mathscr C^\infty$ et l'ellipticit\'e de $dd^c$ impliquent alors qu'il
existe un courant lisse  $\alpha_Z$ sur $\bar G(\C)$, 
(c'est-\`a-dire une forme diff\'erentielle de type $(t,t)$
de classe $\mathscr C^\infty$) et deux courants $u$ et $v$ de types
$(t-1,t)$ et $(t,t-1)$ respectivement de sorte que 
$g_Z=\alpha_Z+\partial u + \bar\partial v$.
Autrement dit, l'\'el\'ement $(0,\alpha_Z)
\in \widehat{\mathrm Z}\vphantom{Z}^{t+1}(\tilde G)$ repr\'esente 
le produit $\tilde D_0\cdots\tilde D_t$.

Soient $\overline{\mathscr F_{1}},\ldots,\overline{\mathscr F_{g}}$
des fibr\'es en droites munis de m\'etriques ad\'eliques sur $\bar G$.
On veut prouver que le produit
\[ \hc_1(\overline{\mathscr L_0})\cdots \hc_1(\overline{\mathscr L_t})\cdot
\hc_1(\overline{\mathscr F_1})\cdots \hc_1(\overline{\mathscr F_g}) = 0 .\]
Par approximation, on peut supposer que les m\'etriques ad\'eliques
sur les $\mathscr F_j$ sont alg\'ebriques donn\'es par des \'el\'ements
$\hc_1(\tilde{\mathscr F_j})\in\hCH^1(\tilde G)$.
Soit $(Z',g_{Z'})$ un repr\'esentant de leur produit.

Si l'on se donne des mod\`eles $\tilde G_n$ de $\bar G$ munis de deux morphismes
$\phi_n:\tilde G_n\ra \tilde G$ et $\psi_n:\tilde G_n\ra\tilde G$
qui \'etendent respectivement l'identit\'e et la multiplication par $n$,
le nombre \`a calculer est la limite des expressions
\begin{align*} P_n & := 
 \frac{1}{n^{t+1}}  \hdeg \psi_n^*(Z,g_Z)\cdot \phi_n^*(Z',g_{Z'}) 
= \frac{1}{n^{t+1}} \hdeg  \psi_n^* (0,\alpha_Z)\cdot \phi_n^*(Z',g_{Z'}) \\
&= \frac{1}{n^{t+1}}  \hdeg  (0,\alpha_Z\circ [n]\wedge \omega_{Z'})
= \frac{1}{n^{t+1}} \int_{\bar G(\C)} (\alpha_Z\circ[n]) \omega_{Z'}  
\end{align*}
d'apr\`es la formule~2.3.1, p.~63 de~\cite{soule-a-b-k92}.

Le fait d'avoir choisi une m\'etrique plate sur les $\mathscr Q_i$
a la cons\'equence suivante. Soit $U\subset A(\C)$ un ouvert
sur lequel les $\mathscr Q_i$ poss\`ede une section de norme~$1$.
Du point de vue diff\'erentiel hermitien,
on peut alors identifier $\pi^{-1}(U)\subset \bar G(\C)$ au produit
$\P^t(\C)\times U(\C)$. 
Sur $\P^t$, le produit d'intersection arithm\'etique analogue $\hat D_0\cdots
\hat D_t$ est repr\'esent\'e par un couple $(0,\alpha)$ o\`u $\alpha$
est une forme diff\'erentielle de type $(t,t)$ sur $\P^t(\C)$.
La forme $\alpha_Z$ s'identifie \`a la forme diff\'erentielle
{\og constante\fg} induite par $\alpha$, et $\alpha_Z\circ[n]$
sur $\bar G(\C)$ s'identifie de m\^eme \`a $[n]^*\alpha$.
La nullit\'e de la hauteur des vari\'et\'es toriques pour les
m\'etriques ad\'eliques canoniques (\cite{maillot97}, proposition 7.1.1)
signifie exactement que
\[ \lim_{n\ra +\infty} n^{-t-1} \int_{\P^t(\C)} [n]^*\alpha = 0 .\]
Par suite, $P_n$ tend vers~$0$ quand $n$ tend vers~$+\infty$,
comme il fallait d\'emontrer.
\end{proof}

\begin{theo}
Relative au fibr\'e inversible m\'etris\'e
$\overline{\mathscr L}\otimes\pi^*\overline{\mathscr M}$, on a
\[ \hat h(\bar G) =
- \frac{1}{(t+1)(t+2) } \left(
       \hat h_{\mathscr M}(q) + \sum_{i=1}^t \hat h_{\mathscr
M}(q-(t+1)q_i) \right). \]
\end{theo}
\begin{proof}
On veut calculer
$\hc_1(\overline{\mathscr L}\otimes\pi^*\overline{\mathscr M})^{g+t+1}$.
L'\'egalit\'e de cycles $[2]_*\overline G=2^{2g+t}\overline G$
et la proposition~\ref{prop:normalisation}
impliquent que si $a\neq t+2$,
\[ \hc_1(\overline{\mathscr L})^a \hc_1(\pi^*\overline{\mathscr
M})^{g+t+1-a} = 0 .\]
Ainsi,
\begin{align*}
\hc_1(\overline{\mathscr L}\otimes\pi^*\overline{\mathscr M})^{g+t+1} 
\hskip -2cm & \hbox{}\hskip 2cm
 = \binom{g+t+1}{t+2} 
  \hc_1(\overline{\mathscr L})^{t+2} \hc_1(\pi^*\overline{\mathscr M})^{g-1}\\
& = \binom{g+t+1}{t+2}
  \big( (t+1)\hc_1 (\overline{\mathscr L_0})
                  + \hc_1(\pi^*\overline{\mathscr Q})\big)^{t+2}
  \hc_1(\pi^*\overline{\mathscr M})^{g-1} \\
& = \binom{g+t+1}{t+2} \sum_{a=0}^{t+2}
       \binom{t+2}a  (t+1)^a \hc_1(\overline{\mathscr L_0})^a
                 \hc_1( \pi^* (\overline{\mathscr Q}))^{t+2-a}
                 \hc_1( \pi^* (\overline{\mathscr M}))^{g-1} \\
&=  \binom{g+t+1}{t+2} \sum_{a=0}^{t+2}
     \binom{t+2}a (t+1)^a
           \pi_* \big( \hc_1(\overline{\mathscr L_0})^a\big)
            \hc_1(\overline{\mathscr Q})^{t+2-a}
            \hc_1(\overline{\mathscr M})^{g-1} 
\end{align*}
d'apr\`es la formule de projection.
On voit que dans cette somme,
si $(t+2-a)+(g-1)>g+1$, c'est-\`a-dire si $a<t$, le terme correspondant
est nul. Il nous faut ainsi calculer
$\pi_* \hc_1(\overline{\mathscr L_0})^a$ pour $a=t$, $a=t+1$ et $a=t+2$.
Ce sont des \'el\'ements de $\hCH^0(A)$, $\hCH^1(A)$ et $\hCH^2(A)$
respectivement.

Le premier est en fait \'egal \`a $\pi_* c_1(\mathscr L_0)^t=1$
(par exemple parce que le produit est repr\'esent\'e num\'eriquement
par le cycle $D_1\cdots D_t$).

Utilisons maintenant la proposition pr\'ec\'edente
que l'on peut d\'evelopper en
\[
 0 =  \hc_1(\overline{\mathscr L_0})^{t+1} + \hc_1(\overline{\mathscr L_0})^t 
        \pi^* \hc_1(\overline{\mathscr Q})
  + \hc_1(\overline{\mathscr L_0})^{t-1}
         \pi^* \left(\sum_{i <j} \hc_1(\overline{\mathscr Q_i})
\hc_1(\overline{\mathscr Q_j})\right)
 + \cdots 
\]
Appliquant $\pi_*$ et utilisant la formule de projection, on obtient ainsi
\[ \pi_* \hc_1(\overline{\mathscr L_0})^{t+1}
      = - \hc_1(\overline{\mathscr Q}) \]
tandis que si on multiplie par $\hc_1(\overline{\mathscr L_0})$
avant d'appliquer $\pi_*$, on trouve
\[
 \pi_* \hc_1(\overline{\mathscr L_0})^{t+2}
  = \hc_1(\overline{\mathscr Q}) ^2
           - \sum_{i<j} \hc_1(\overline{\mathscr Q_i})
             \hc_1(\overline{\mathscr Q_j}) 
= \frac12 \hc_1(\overline{\mathscr Q})^2
 +\frac 12\sum_{i=1}^t \hc_1(\overline{\mathscr Q_i})^2.
\]
(\`A proprement parler, on n'\'ecrit ces \'egalit\'es qu'apr\`es
avoir multipli\'e les deux membres par $\hc_1(\overline{\mathscr
Q})^a\hc_1(\overline{\mathscr M})^{b}$ pour
$a+b=g$ ou $g-1$.)
Il en r\'esulte que
\begin{align*}
\hc_1(\overline{\mathscr L}\otimes\pi^*\overline{\mathscr M})^{g+t+1}
\hskip-2cm \\
&= \binom{g+t+1}{t+2} \left(
\frac{(t+2)(t+1)}2 (t+1)^t \hc_1(\overline{\mathscr Q})^2 
- (t+2) (t+1)^{t+1} \hc_1(\overline{\mathscr Q})^2
\vphantom{\sum_i}\right. \\
& \qquad\qquad\left.
   {}+ (t+1)^{t+2} \frac12 \left(
   \hc_1(\overline{\mathscr Q})^2 +
   \sum_{i} \hc_1(\overline{\mathscr Q_i})^2
        \right)\right)
 \hc_1(\overline{\mathscr M})^{g-1} \\
&= \frac12 \binom{g+t+1}{t+2} (t+1)^{t+1} 
   \left( - \hc_1(\overline{\mathscr Q})^2 
         + (t+1) \sum_{i=1}^t \hc_1(\overline{\mathscr Q_i})^2 \right)
   \hc_1(\overline{\mathscr M})^{g-1} \\
&= - \frac1g \binom{g+t+1}{t+2} (t+1)^{t+1}
        c_1(\mathscr M)^g
    \left(
           (t+1) \sum_{i=1}^t \hat h_{\mathscr M}(q_i)
           - \hat h_{\mathscr M}\big(\sum_{i=1}^t q_i\big) \right)
\end{align*}
d'apr\`es le th\'eor\`eme~\ref{theo:Q2L}.
Ainsi, on a
\begin{multline*}
 \hc_1(\overline{\mathscr L}\otimes\pi^*\overline{\mathscr M})^{g+t+1} = \\
{} =  
      - \frac{(g+t+1)!}{g!(t+2)!} (t+1)^{t+1}
        c_1(\mathscr M)^g
    \left(
           (t+1) \sum_{i=1}^t \hat h_{\mathscr M}(q_i)
           - \hat h_{\mathscr M}\big(\sum_{i=1}^t q_i\big) \right)
\end{multline*}

Le calcul de l'intersection g\'eom\'etrique est plus facile.
Comme les $\mathscr Q_i$ sont alg\'ebriquement \'equivalent \`a z\'ero,
tout se passe du point de vue num\'erique comme si l'on avait un produit
et
\[ c_1(\mathscr L\otimes\pi^*\mathscr M)^{g+t}
        = \binom{g+t}t (t+1)^t c_1(\mathscr M)^g. \]

Finalement,
\begin{align*}
\hat h_{\overline{\mathscr L}\otimes\pi^*\overline{\mathscr M}}(\bar G)
 & = \frac{1}{g+t+1}
     \frac{
        \hc_1(\overline{\mathscr L}\otimes\pi^*\overline{\mathscr M})^{g+t+1}}
       {c_1(\mathscr L\otimes\pi^*\mathscr M)^{g+t}} \\
& = - \frac{1}{g+t+1} \frac{(g+t+1)!}{g!(t+2)!} \frac{g! t!}{(g+t)!}
         \frac{(t+1)^{t+1} c_1(\mathscr M)^g}{(t+1)^t c_1(\mathscr M)^g}
\times\\
& \qquad\qquad\times
      \left( 
           (t+1) \sum_{i=1}^t \hat h_{\mathscr M}(q_i)
           - \hat h_{\mathscr M}\big(\sum_{i=1}^t q_i\big) \right) \\
&= - \frac{1}{t+2} \left(  (t+1) \sum_{i=1}^t \hat h_{\mathscr M}(q_i)
           - \hat h_{\mathscr M}\big(\sum_{i=1}^t q_i\big) \right).
\end{align*}
En utilisant la quadraticit\'e de la hauteur de N\'eron--Tate,
on a enfin
\begin{multline*}
\hat h(q)+\sum_{i=1}^t \hat h (q-(t+1)q_i)
= (t+1)\hat h(q) - 2(t+1) \hat h(q) + (t+1)^2 \sum_{i=1}^t \hat h(q_i)\\
= (t+1) \left( (t+1)\sum_{i=1}^t \hat h(q_i) - \hat h(q) \right).
\end{multline*}
si bien que le th\'eor\`eme est d\'emontr\'e.
\end{proof}

Le corollaire que voici acquerra tout son sel
au paragraphe suivant.

\begin{coro}
On a $\hat h(\bar G)\leq 0$, avec \'egalit\'e si et seulement
si la vari\'et\'e semi-ab\'elienne $G$ est isotriviale.
\end{coro}
\begin{proof}
Comme la hauteur de N\'eron--Tate est d\'efinie positive sur $A(k)\otimes\R$,
la formule du th\'eor\`eme pr\'ec\'edent montre que $\hat h(\bar G)\leq 0$.
De plus, on a \'egalit\'e si et seulement si $\hat h(q)=0$ et
$\hat h(q-(t+1)q_i)=0$ pour tout $i$, 
c'est-\`a-dire, $0=q=(t+1)q_i$ dans $A(k)\otimes\R$.
Par suite tous les $q_i$ sont de torsion dans $A(k)$,
les fibr\'es en droites $\mathscr Q_i$ d'ordre fini et
$G$ est isotriviale.
\end{proof}

\begin{rema}
On peut effectuer le calcul de la hauteur normalis\'ee pour le fibr\'e inversible
\[ \mathscr L' = \bigotimes_{0\leq i\leq t}
 \mathscr L_i^{n_i} \otimes \pi^*\mathscr M \]
et, posant 
\[ r= \big( n_1q_1+\cdots+n_t q_t\big) 
      -\big( n_0+\cdots+n_t\big) \frac{ q_1+\cdots+q_t }{t+1} \]
dans $ A(k)\otimes\R$,
on obtient la formule
\[ \hat h'(\bar G) = \left(\frac{n_0+\cdots+n_t}{t+1}\right)^2 \hat h(\bar G)
        - \hat h_{\mathscr M} (r). \]
Ainsi, parmi les fibr\'es en droites $G$-lin\'earis\'es m\'etris\'es
$\overline{\mathscr L}$ sur $\bar G$ tels que $\mathscr
L\otimes\pi^*\mathscr M$ est de degr\'e donn\'e, c'est le fibr\'e anticanonique
relatif qui donne \`a $\bar G$ la plus grande hauteur.
\end{rema}

\begin{lemm}
On a 
\[ \mu_{\overline{\mathscr L}}(\bar G)=
 - \max \big( \hat h_{\mathscr M}(q), \max_i \big( \hat h_{\mathscr
   M}(q-(t+1)q_i)\big)
   \big). \]
\end{lemm}
\begin{proof}
L'intersection des $D_i$ pour $0\leq i\leq t$ est vide.
Par suite, il suffit de prouver les \'egalit\'es
\[ \min_{x\not\in D_0}
   \hat h_{\overline{\mathscr L}\otimes\pi^*\overline{\mathscr M}}(x)
 = -\hat h_{\overline{\mathscr M}}(q), \quad
   \min_{x\not\in D_i}
   \hat h_{\overline{\mathscr L}\otimes\pi^*\overline{\mathscr M}}(x)
 = -\hat h_{\overline{\mathscr M}}(q-(t+1)q_i). \]
Compte-tenu de l'isom\'etrie
$\mathscr L_i \simeq \mathscr L_0\otimes\pi^*\mathscr Q_i$,
on a pour $x\in\bar G(\bar\Q)$,
\begin{align*}
 \hat h_{\overline{\mathscr L}\otimes\pi^*\overline{\mathscr M}}(x)
    & = (t+1) h_{\overline{\mathscr L_0}} (x)
    + h_{\overline{\mathscr Q}}(\pi(x)) + h_{\overline{\mathscr M}}(\pi(x)) \\
   & = (t+1) h_{\overline{\mathscr L_0}}(x)
  + h_{\overline{\mathscr M}}(\pi(x)+q) - h_{\overline{\mathscr M}}(q)
\end{align*}
en vertu de l'isomorphisme de faisceaux inversibles m\'etris\'es
du corollaire~\ref{coro:Q-L}.
Si $x$ n'appartient pas au support de $D_0$, l'argument de la preuve
du lemme~\ref{lemm:positivite} implique que
$h_{\overline{\mathscr L_0}}(x)\geq 0$.
Ainsi, pour $x\not\in D_0$, on a
\[ h_{\overline{\mathscr L}\otimes\pi^*\overline{\mathscr M}} (x)
            \geq - h_{\overline{\mathscr M}}(q). \]
L'intersection des $D_i$ pour $i\neq 0$ est un sous-sch\'ema ferm\'e
de $\bar G$ que la projection $\pi$ identifie canoniquement \`a $A$
et la restriction de $\mathscr L_0$ s'identifie alors au fibr\'e
trivial. Ainsi, si $x$ est dans tous les $D_i$ mais
pas dans $D_0$, et si de plus $\pi(x)=-q$, on a 
$  h_{\overline{\mathscr L}\otimes\pi^*\overline{\mathscr M}} (x)
= - h_{\overline{\mathscr M}}(q)$.

De la m\^eme mani\`ere, on prouve
que si $x\not\in D_i$,
\[ \hat h_{\overline{\mathscr L}\otimes\pi^*\overline{\mathscr M}} (x)
            \geq - \hat h_{\overline{\mathscr M}}(q-(t+1)q_i) \]
et que cette minoration est optimale sur $\bar G\setminus D_i$.
\end{proof}

\begin{coro}
On a $\mu(\bar G)\leq 0$, avec \'egalit\'e si et seulement
si $G$ est une vari\'et\'e semi-ab\'elienne isotriviale.
\end{coro}

\section{Variantes de la conjecture de Bogomolov}
\label{sect:conj}

Soient $G$ une vari\'et\'e semi-ab\'elienne sur un corps
de nombres $k\subset\bar\Q$. On fixe une compactification $\bar G$
ainsi qu'une hauteur normalis\'ee comme pr\'ec\'edemment.
Soit $X$ une sous-vari\'et\'e irr\'eductible de $G$ et
$\bar X$ son adh\'erence dans $\bar G$.
Rappelons qu'on appelle sous-vari\'et\'e de torsion de
$G$ toute sous-vari\'et\'e de $G\otimes\bar \Q$
qui est la translat\'ee d'un sous-groupe alg\'ebrique
de $G\otimes\bar \Q$ par un point de torsion.

F.~Bogomolov \'enonce dans~\cite{bogomolov80b}
une conjecture sur la hauteur de N\'eron--Tate des points
alg\'ebrique des sous-vari\'et\'es d'une vari\'et\'e ab\'elienne.
Cet \'enonc\'e se g\'en\'eralise aux vari\'et\'es
semi-ab\'eliennes de la fa\c{c}on suivante.

\begin{enonce}{\'Enonc\'e}%
   [G\'en\'eralisation de la conjecture de Bogomolov]
   \label{conj:bogo}
\begin{enumerate}
\item L'ensemble des sous-vari\'et\'es de torsion contenues dans $X$
n'a qu'un nombre fini d'\'el\'ements maximaux.
\item Si $X^*$ d\'esigne le compl\'ementaire dans $X$ des sous-vari\'et\'es
de torsion contenues dans $X$, il existe un r\'eel $\eps>0$ tel que
pour tout $x\in X^*(\bar\Q)$, $\hat h(x)>0$.
\end{enumerate}
\end{enonce}

Ainsi que l'a montr\'e Zhang dans~\cite{zhang95},
cet \'enonc\'e \'equivaut
dans le cas des vari\'et\'es ab\'eliennes
\`a un autre \'enonc\'e que Philippon avait
conjectur\'e ind\'ependamment~\cite{philippon91}.

\begin{enonce}{\'Enonc\'e}%
   [G\'en\'eralisation de la conjecture de Philippon]
   \label{conj:pph}
La vari\'et\'e $\bar X$ est de hauteur normalis\'ee nulle si et seulement si
$X$ est une sous-vari\'et\'e de torsion de~$G$.
\end{enonce}

La conjecture d'\'equir\'epartition que nous \'enon\c{c}ons ci-dessous
a \'et\'e introduite par Szpiro, Ullmo et Zhang dans~\cite{szpiro-u-z97}
lorsque $G$ est une vari\'et\'e ab\'elienne ; il y est prouv\'e
qu'elle \'equivaut dans ce cas \`a la conjecture
de Bogomolov pr\'ec\'edente. L'\'equir\'epartition des suites de petits points
joue un r\^ole essentiel dans les d\'emonstrations par Ullmo et Zhang 
des analogues ab\'eliens des \'enonc\'es~\ref{conj:bogo} et~\ref{conj:pph}.

\begin{enonce}{\'Enonc\'e}%
  [Conjecture d'\'equir\'epartition des petits points]
  \label{conj:equirep} 
Soit $(x_n)$ une suite de points de $G(\bar \Q)$ 
dont la hauteur normalis\'ee tend vers~$0$ et telle qu'aucune sous-suite
ne soit contenue dans un sous-groupe alg\'ebrique strict de~$G$.
D\'esignons
par $O(x_n)$ l'adh\'erence sch\'ematique de $x_n$ dans $G$
et $\delta_{O(x_n)}$ la mesure de comptage sur l'ensemble fini $O(x_n)(\C)$.
Alors, la suite de mesures de probabilit\'e
$ \frac{1}{\# O(x_n)} \delta_{O(x_n)(\C)}$
converge vaguement
vers la mesure de Haar normalis\'ee port\'ee par le sous-groupe
compact maximal de $G(\C)$.
\end{enonce}

Lorsque $X$ est une vari\'et\'e ab\'elienne, ces trois \'enonc\'es
sont vrais (voir~\cite{ullmo98}, \cite{zhang98} et~\cite{david-p98}).
Lorsque $X$ est un tore, les deux premiers \'enonc\'es ont \'et\'e d\'emontr\'e
par Zhang dans~\cite{zhang95b} puis par d'autres m\'ethodes
par
Bombieri--Zannier~\cite{bombieri-z95} et
Schmidt~\cite{schmidt96}.
L'\'enonc\'e d'\'equir\'epartition est dans ce cas un th\'eor\`eme de
Y.~Bilu~\cite{bilu97}.

Dans le cas semi-ab\'elien, ces \'enonc\'es ne sont pas forc\'ement
\'equivalents et nous allons m\^eme donner un contre-exemple
\`a l'\'enonc\'e~\ref{conj:pph} un peu plus bas.
N\'eanmoins, si $\mu(\bar X)\geq 0$, c'est-\`a-dire si la hauteur
de tous les points de $\bar G$ est positive, on peut dire quelque chose.

\begin{lemm}
Si $\mu(\bar X)\geq 0$, les \'enonc\'es~\ref{conj:bogo} et~\ref{conj:pph}
sont \'equivalents.
\end{lemm}
\begin{proof}
Bien qu'il faille juste reprendre les arguments utilis\'es dans l'\'equivalence
des analogues ab\'eliens de ces \'enonc\'es,
nous d\'etaillons la preuve pour le confort du lecteur.

Lorsque $\mu(\bar X)\geq 0$, l'in\'egalit\'e du th\'eor\`eme~\ref{theo:1.10}
s'\'ecrit en effet
\[ e(X) \geq \hat h(\bar X) \geq \frac{1}{1+\dim X} e(X) \]
si bien que la nullit\'e de $\hat h(\bar X)$ ou celle de $e(X)$
sont \'equivalentes.
Le reste de la preuve proc\`ede comme dans~\cite{zhang95b} :

Supposons que l'\'enonc\'e~\ref{conj:bogo} est v\'erifi\'e.
Si $X$ n'est pas une sous-vari\'et\'e de torsion, $X^*$ est non 
vide si bien que $e(X)\geq \eps>0$ et donc $\hat h(X)>0$.
R\'eciproquement, si $X$ est une sous-vari\'et\'e de torsion, 
ses points de torsion sont Zariski-dense, ce qui implique
$e(X)=0$ et dont $\hat h(\bar X)=0$.

Supposons maintenant l'\'enonc\'e~\ref{conj:pph} satisfait.
Si~$X$ n'est pas une sous-vari\'et\'e de torsion,
donc si $\hat h(\bar X)>0$,
il existe un ouvert de Zariski non vide $U$ de $X$ tel que si $x\in U$,
la hauteur de $x$ est strictement positive. En particulier, $U$
ne rencontre aucune sous-vari\'et\'e de torsion de $G$ et toute sous-vari\'et\'e
de torsion de $G$ contenue dans $X$ est contenue dans $X\setminus U$.
Par r\'ecurrence descendante sur la dimension des composantes irr\'eductibles
de $X\setminus U$, on voit que les sous-vari\'et\'es de torsion maximales
contenues dans $X$ sont en nombre fini.
Si de plus la borne inf\'erieure de la hauteur normalis\'ee
sur l'ouvert compl\'ementaire $X^*$ \'etait nulle, on pourrait
en d\'eduire l'existence d'une sous-vari\'et\'e de hauteur nulle qui n'est
pas incluse dans les pr\'ec\'edentes, et donc d'apr\`es
l'\'enonc\'e~\ref{conj:pph}, une autre sous-vari\'et\'e de torsion, ce qui 
est une contradiction.
\end{proof}

\begin{lemm}
Si $\mu(\bar X)\geq 0$, l'\'enonc\'e~\ref{conj:equirep}
implique l'\'enonc\'e~\ref{conj:pph}.
\end{lemm}
\begin{proof}
Supposons en effet que la hauteur normalis\'ee de $\bar X$ est nulle
mais que $X$ n'est pas une sous-vari\'et\'e de torsion.
Comme on a suppos\'e que $\mu(\bar X)\geq 0$,
le th\'eor\`eme~\ref{theo:1.10} implique comme pr\'ec\'edemment
que $e(X)=0$, c'est-\`a-dire
qu'il existe une suite $(x_n)$ de points de
$X(\bar\Q)$ dont la hauteur tend vers $0$ et convergeant vers
le point g\'en\'erique de $X$. En particulier, aucune sous-suite n'est
contenue dans un sous-groupe alg\'ebrique strict.
L'\'enonc\'e~\ref{conj:equirep} implique alors que la suite
des mesures $\delta_{O(x_n)(\C)}/\#O(x_n)$ converge vaguement
vers la mesure de Haar normalis\'ee port\'ee
par le sous-groupe compact maximal de $G(\C)$. Par suite, $X(\C)$
contient ce sous-groupe compact maximal. 
Remarquons alors que le sous-groupe compact maximal de $\gm(\C)$
est Zariski-dense dans $\gm$ ;
par suite, le sous-groupe compact maximal d'un tore complexe
est aussi Zariski-dense, ce qui implique que le sous-groupe
compact maximal de $G(\C)$ est dense dans chaque fibre
de $G\ra A$. Puisque $X$ est une sous-vari\'et\'e alg\'ebrique,
elle doit contenir toutes ces fibres, et donc \^etre \'egal \`a $\bar G$,
d'o\`u une contradiction.
\end{proof}

\bigskip

Il r\'esulte des calculs du paragraphe~\ref{sect:calcul}
que les vari\'et\'es semi-ab\'eliennes qui ne sont pas isotriviales
sont de hauteur normalis\'ee strictement n\'egative. 
On a autrement dit le th\'eor\`eme suivant :
\begin{theo} \label{theo:faux}
La g\'en\'eralisation~\ref{conj:pph} de la conjecture de Philippon 
est \emph{fausse}.
\end{theo}

N\'eanmoins, et m\^eme lorsque $G$ n'est pas isotriviale, 
$\mu(\bar G)<0$ et il est donc possible (probable ?)
que les deux autres conjectures soient v\'erifi\'ees.

Dans la suite de ce texte, nous allons \emph{d\'emontrer}
ces trois conjectures lorsque $G$ est isotriviale.
Dans ce cas, $\mu(\bar G)=0$ si bien que l'on a 
d'ores et d\'ej\`a les implications
\[ \text{Conj. de Bogomolov} \Leftrightarrow \text{Conj. de Philippon}
    \Leftarrow \text{Conj. d'\'equir\'epartition}. \]
De plus, il est clair que si les conjectures~\ref{conj:bogo}
et~\ref{conj:equirep} sont vraies pour une vari\'et\'e semi-ab\'elienne
$G$, elles sont vraies pour toute vari\'et\'e semi-ab\'elienne qui
lui est isog\`ene.
Par suite, il suffit de traiter le cas d'une vari\'et\'e semi-ab\'elienne
qui est un produit $A\times T$. La th\'eorie du \S\,\ref{sect:constr}
se simplifie grandement dans ce cas et on a $\bar G=A\times X$.

\section{Un th\'eor\`eme d'\'equir\'epartition g\'en\'erique
en rang torique \'egal \`a un}
\label{sect:generique}

Le th\'eor\`eme principal de cette section est le
th\'eor\`eme~\ref{theo:equirep} ci-dessous. Dans le cas d'un tore,
il s'agit du th\'eor\`eme de Bilu, dont nous donnons en appendice une
preuve arakelovienne.
Nous g\'en\'eralisons en fait cette preuve
et  obtenons ainsi un th\'eor\`eme d'\'equir\'epartition g\'en\'erique
pour les vari\'et\'es semi-ab\'eliennes de rang torique \'egal
\`a~$1$. Une adaptation des m\'ethodes employ\'ees par Ullmo
puis Zhang pour en d\'eduire la conjecture de Bogomolov pour
les vari\'et\'es ab\'eliennes, ainsi que de la m\'ethode employ\'ee par Bilu
pour passer de $\gm$ \`a un tore quelconque nous permettra
de prouver la conjecture de Bogomolov pour les vari\'et\'es semi-ab\'eliennes
isotriviales au paragraphe suivant.

\begin{theo}\label{theo:equirep}
Soient, sur un corps de nombres $k$, une vari\'et\'e ab\'elienne
$A$ et $G$ la vari\'et\'e semi-ab\'elienne $A\times\gm$ de
compactification $\bar G=A\times\P^1$. Notons $\pi:\bar G\ra A$
la projection naturelle.
Soit $\overline{\mathcal L}$ 
un fibr\'e ample sur $\bar G$, muni de sa m\'etrique ad\'elique
canonique qui d\'efinit une hauteur normalis\'ee sur $G$ au sens de la
d\'efinition~\ref{defi:h-norm}.

Soit $V\subset \bar G$ une sous-vari\'et\'e ferm\'ee
(irr\'eductible et de dimension $>0$) qui rencontre $G$.
Soit $(x_n)$ une suite {\em g\'en\'erique\/} de {\em petits points\/}
de $V(\bar k)$, c'est-\`a-dire qu'on suppose :
\begin{itemize}
\item qu'aucune sous-suite de $(x_n)$ n'est contenue dans un ferm\'e
strict de $V$ ;
\item et que $h_{\overline{\mathcal L}}(x_n)\ra 0$ quand $n\ra+\infty$.
\end{itemize}
On suppose de plus que l'une des deux hypoth\`eses suivantes est
satisfaite :
\begin{itemize}
\item ou bien $V=\bar G$ ;
\item ou bien la projection naturelle $V(\C)\ra A(\C)$
est g\'en\'eriquement finie.
\end{itemize}

Fixons une place $\sigma:k\hra\C $. Pour tout $n$, soit $O(x_n)$
l'orbite de $x_n$ dans $\bar G$ sous l'action de $\Gal(\bar k/k)$
et $\mu_{O(x_n)}$ la mesure
de probabilit\'e sur $V(\C_\sigma)$ port\'ee par $O(x_n)$ et invariante
sous l'action de $\Gal(\C_\sigma /k)$. 

{\em Alors,} la suite des mesures $\mu_{O(x_n)}$ converge vaguement vers
la mesure de probabilit\'e $\nu_{\overline{\mathcal L}|V}$,
d\'efinie dans la d\'efinition~\ref{defi:bt}.
\end{theo}

Avant de d\'emontrer le th\'eor\`eme, remarquons que l'on
peut supposer
\[ \overline{\mathscr L}
   =\overline{\mathscr O_{\P}(1)}\otimes\pi^*\overline{\mathscr M}. \]
On notera $D$ le diviseur {\og infini\fg} $A\times\{\infty\}$ de $\bar G$.
On a $\mathcal O_{\bar G}(D)=\mathcal O_{\P}(1)$, d'o\`u une section canonique
$\mathsf s_D$ de $\mathcal O_{\P}(1)$.

Le r\'esultat clef dans la d\'emonstration de ce genre
de th\'eor\`eme d'\'equir\'epartion par des m\'ethodes de g\'eom\'etrie
d'Arakelov est le suivant, d\'emontr\'e par Szpiro, Ullmo
et Zhang quand les m\'etriques hermitiennes sont lisses : 
\begin{prop}[Szpiro, Ullmo, Zhang] \label{prop:suz}
Consid\'erons un fibr\'e en droites $\overline{\mathcal L}$ ample, muni
d'une m\'etrique ad\'elique ample
sur une vari\'et\'e projective $X$ d\'efinie
sur un corps de nombres $k\subset\C$.
On suppose que $X$ est de hauteur nulle pour $\overline{\mathcal L}$.
Soient $(x_n)$ une suite g\'en\'erique de petits points dans $X(\bar\Q)$
et $f:X(\C)\ra\R$ une fonction continue telle que la propri\'et\'e
suivante soit satisfaite :
\[  \text{pour tout $\eps>0$ assez petit, le courant
    $ c_1(\overline{\mathcal L})+ \eps dd^c f $
est positif sur $X$.}  \]
Alors, on a l'in\'egalit\'e
\[  \liminf_{n\ra\infty} \int_X f \mu_{O(x_n)} \geq \int_X f
\nu_{\overline{\mathcal L}}.  \]
\end{prop}
\begin{proof}
Pour $\eps$ assez petit, comme dans l'\'enonc\'e de la proposition,
consid\'erons le fibr\'e en droites m\'etris\'e $\overline{\mathcal L}(\eps f)$
sur $X$ qui n'est en fait autre que $\overline{\mathcal L}$,
mais dont la m\'etrique hermitienne \`a l'infini ($k\subset\C$)
est multipli\'ee par $\exp(-\eps f)$.
Par d\'efinition, le courant de courbure de $\overline{\mathcal L}(\eps f)$ est
\'egal \`a $c_1(\overline{\mathcal L})+\eps dd^c f$, donc positif ;
en particulier, la m\'etrique ad\'elique
de $\overline{\mathcal L}(\eps f)$ est semi-positive,
si bien qu'on peut
appliquer le th\'eor\`eme~\ref{theo:1.10} \`a la suite
g\'en\'erique $(x_n)$ et au fibr\'e m\'etris\'e $\overline{\mathcal L}(\eps f)$.
On obtient ainsi
\[  \liminf_{n\ra\infty} h_{ \overline{\mathcal L}(\eps f) }(x_n)
        \geq h_{\overline{\mathcal L}(\eps f)}(X).  \]
Calculons alors les deux membres de cette in\'egalit\'e.
Pour le membre de gauche, on a :
\[  h_{ \overline{\mathcal L}(\eps f) }(x_n)
  = h_{\overline{\mathcal L}} (x_n) +\eps \int_X f\mu_{O(x_n)},  \]
tandis que le membre de droite vaut d'apr\`es la
proposition~\ref{prop:calcul}
\[  h_{\overline{\mathcal L}} (X) + \eps \int_X f \nu_{\overline{\mathcal L}}
 + O(\eps^2).  \]

Comme $X$ est de hauteur nulle pour $\overline{\mathcal L}$
et comme  $(x_n)$ est une suite de petits points, il en r\'esulte
en divisant par $\eps$
l'in\'egalit\'e
\[  \liminf_{n\ra\infty} \int_X f \mu_{O(x_n)} \geq \int_X f
\nu_{\overline{\mathcal L}} + O(\eps).  \]
Il reste \`a faire tendre $\eps$ vers~$0$ par valeurs sup\'erieures
et la proposition est d\'emontr\'ee.
\end{proof}

On rassemble en un lemme des r\'esultats connus 
(voir par exemple Zhang, \cite{zhang95}, 6.4).
\begin{lemm}\label{lemm:norm}
La norme de la section canonique $\mathsf s_D$ de $\mathcal O_{\P}(1)$ est
\[  \norm{\mathsf s_D}(z,[t:u])= \frac{|u|}{\max(|t|,|u|)},
            \quad z\in A(\C), \quad [t:u]\in\P^1(\C),  \]
et, notant $\zeta=t/u\in\C\union\{\infty\}$
on a pour le courant de courbure la formule :
\[  c_1(\overline{\mathcal O_{\P}(1)}) =
\frac{d\Arg(\zeta)}{2\pi}\wedge\delta_{|\zeta|=1}.  \]
\end{lemm}
On a not\'e $[t:u]$ l'unique point de $\P^1(\C)$ 
dont les coordonn\'ees homog\`enes sont $(t,u)$. De plus, 
si $X$ est une vari\'et\'e analytique complexe, $i:Y\hra X$
une sous-vari\'et\'e diff\'erentielle r\'eelle de $X$
et $\eta$ une forme diff\'erentielle lisse dans un voisinage de $Y$,
on note $\eta\wedge\delta_{Y}$ le {\em courant} (port\'e par $Y$)
$i_*i^*\eta$.

\begin{proof}
Comme nous sommes dans une situation produit,
le calcul se ram\`ene au calcul analogue sur $\P^1$, avec $D=[1:0]=\infty$.
La m\'etrique de Fubini-Study sur $\mathscr O(D)$ donne \`a $\mathsf s_D$ la norme
\[ \norm{\mathsf s_D}_{FS} ([t:u]) = \frac{\abs{u}}{(\abs t^2+\abs u^2)^{1/2}}. \]
L'isomorphisme $[n^*]\mathscr O(D)\simeq \mathscr O(D)^n$ envoie
la section $[n]^*\mathsf s_D$ sur la section $\mathsf s_D^n$.
Le proc\'ed\'e limite par lequel $\norm{\mathsf s_D}$ est d\'efinie implique alors que
\[ \norm{\mathsf s_D}([t:u]) = \lim_{n\ra\infty} \norm{\mathsf s_D}_{FS}([t^n:u^n])
^{1/n}
  = \lim_n \frac{\abs u}{(\abs t^{2n}+\abs u^{2n})^{1/2n})}
  = \frac{\abs u}{\max(\abs u,\abs t)}, \]
comme annonc\'e.

D'autre part, $c_1(\overline{\mathscr O(D)})$ est 
la mesure $dd^c \log\norm{\mathsf s_D}+\delta_\infty$ sur $\P^1(\C)$.
Comme $\log\norm{\mathsf s_D}$ est, hors du sous-groupe compact maximal $S^1$,
le logarithme de la valeur absolue d'une fonction m\'eromorphe,
cette mesure est support\'ee par le cercle. Elle est de plus invariante
par rotation et de masse totale~$1$. 
Par suite, l'unicit\'e de la mesure de Haar implique
que 
\[ c_1(\overline{\mathscr O(D)})=\frac{1}{2\pi}d\Arg(\zeta)\wedge
 \delta_{|\zeta|=1}.\]
Le lemme est donc d\'emontr\'e. 
\end{proof}

\begin{lemm}\label{lemm:Phi}
Soit $f$ une fonction $\mathcal C^\infty$ sur $\bar G(\C)$. Il existe
alors une fonction $\Phi$ sur $\bar G(\C)$ ayant les propri\'et\'es
suivantes :
\begin{itemize}
\item $\Phi$ est continue ;
\item $\Phi=f$ sur le sous-groupe compact maximal de $G(\C)$ ;
\item pour tout $\eps>0$ assez petit, $c_1(\overline{\mathcal L})+\eps dd^c\Phi$
est un courant de type $(1,1)$ {\em positif\/} sur~$\bar G$.
\end{itemize}
\end{lemm}
\begin{proof}
Avec les notations du lemme~\ref{lemm:norm},
il nous suffit de construire une fonction continue $\Phi$ sur $A\times
\P^1$ (dont les coordonn\'ees homog\`enes sont toujours not\'ees $[t:u]$)
v\'erifiant les propri\'et\'es suivantes :
\begin{itemize}
\item $\Phi$ est continue ;
\item pour tout $z\in A$ et tout $\theta\in\R$, 
$\Phi(z,[e^{i\theta}:1])=f(z,[e^{i\theta}:1])$ ;
\item pour tout $\eps$ assez petit, $c_1(\overline{\mathcal L})+\eps
dd^c\Phi$ est un courant de type $(1,1)$ positif.
\end{itemize}

Posons, pour $\lambda\in\R$,
\[  \Phi(z,[t:u])= f(z,[t:u]) + \lambda \log
\frac{|t|^2+|u|^2}{2\max(|t|^2,|u|^2)},
\] 
de sorte que les deux premi\`eres conditions sont satisfaites.
On a, notant $\zeta=t/u=r e^{i\theta}$,
\[  dd^c\Phi = dd^c f 
+ \lambda \frac{i}{2\pi (1+r^2)^2} d\zeta\wedge d\bar \zeta
- \lambda \frac{d\theta}{2\pi}\wedge \delta_{r=1}.  \]

D'autre part, 
on a une expression de la forme
\[  c_1(\overline{\mathcal L})= \pi^* \omega_A
 + \frac{d\theta}{2\pi} \wedge \delta_{r=1} ,  \]
o\`u $\omega_A=c_1(\overline{\mathscr M})$
est une forme diff\'erentielle de type $(1,1)$
lisse sur $A$ et strictement positive.
Ainsi,
\[  c_1(\overline{\mathcal L})+\eps dd^c\Phi
= \left(\pi^*\omega_A + \eps\lambda\frac{i}{2\pi(1+|t|^2)^2} d\zeta\wedge
d\bar \zeta
 + \eps dd^c f\right)
 + (1-\eps\lambda) \frac{d\theta}{2\pi} \wedge \delta_{r=1}.  \]
Nous allons montrer qu'il existe $\lambda_0>0$ et $\eps_0>0$
tel que si $\eps<\eps_0$ et $\lambda=\lambda_0$, le premier
terme est une forme diff\'erentielle de type $(1,1)$, lisse
et strictement positive. Si de plus $\eps<1/\lambda_0$, le second
terme est un courant positif, ce qui conclura la d\'emonstration
du lemme.

Le premier terme s'interpr\`ete en effet
comme une forme hermitienne sur le fibr\'e
tangent de $\bar G$. La restriction de ce fibr\'e aux ouverts
$U_1$ et $U_2$ de $\bar G$ d\'efinis par $|t|<2$
et $|t|>1/2$ est triviale, une base du fibr\'e dual \'etant donn\'ee
par l'adjonction de $dt$ (resp.\ de $d(1/t)$) \`a une base de $\Omega^1_A$.
Dans ces ouverts, cette forme hermitienne s'\'ecrit ainsi par blocs
\[  \begin{pmatrix}
        h_A(z) + \eps a(z,t)  & \eps b(z,t)  \\
        \eps b^*(z,t) & \eps c(z,t)  + \eps\lambda  \kappa(t) 
    \end{pmatrix},  \]
o\`u $a$, $b$ et $c$ se d\'eduisent des coordonn\'ees de $dd^c f$
dans les bases locales du fibr\'e cotangent, tandis que
$\kappa(t)$ vaut $\big(2\pi(1+|t|^2)^2\big)^{-1}$ dans $U_1$
et $\big(2\pi(1+|1/t|^2)^2 \big)^{-1}$ dans $U_2$. Les coefficients
de cette matrice sont en particulier des fonctions {\em continues\/}
de $(z,t)\in U_1$ (resp.~$U_2$).

Une telle matrice hermitienne est d\'efinie positive si
et seulement si tous ses mineurs principaux sont positifs.
Par relative compacit\'e des deux ouverts $U_1$ et $U_2$,
il existe ainsi une constante $\eps_1>0$ telle que pour
$|\eps|<\eps_1$, les $g$ premiers mineurs sont $>0$ en tout
point de $\bar G$. Si $|\eps|<\eps_1$, il reste \`a s'assurer 
de la positivit\'e du discriminant de la forme hermitienne.
Or, le discriminant s'\'ecrit, en d\'eveloppant suivant la derni\`ere colonne,
\[  \eps \big( \lambda \kappa(t) + c(z,t) \big)  \det(h_A+\eps a(z,t) ) 
 + O (\eps^2),  \]
le $O$ \'etant uniforme sur $U_1$ et $U_2$.
Par compacit\'e, il existe $\lambda_0$ tel que $\lambda_0 \kappa(t)+c(z,t)>0$
en tout point de $\bar G$. Alors, lorsque $\eps>0$ tend vers $0$,
le premier terme domine et le discriminant est bien strictement
positif d\`es que $\eps>0$ est assez petit.
\end{proof}

\begin{proof}[Preuve du th\'eor\`eme~\ref{theo:equirep} quand $V=\bar G$]
Soit $f:\bar G\ra\R$ une fonction continue que le th\'eor\`eme
de Stone-Weierstra\ss\ nous permet en fait de supposer 
$\mathcal C^\infty $ sur $\bar G$.
Choisissons une fonction $\Phi$ comme dans le lemme~\ref{lemm:Phi}.
En lui appliquant la proposition~\ref{prop:suz},
il en r\'esulte l'in\'egalit\'e
\begin{equation}\tag{$\star$}\label{ineq}
\liminf_n \int_{\bar G}\Phi \mu_{O(x_n)}
\geq \int_{\bar G} \Phi \nu_{\overline{\mathcal L}}. 
\end{equation}

\begin{lemm}
La mesure $c_1(\overline{\mathcal L})^{\dim G}$ (d\'efinie par
la proposition--d\'efinition~\ref{defi:bt})
est un multiple de la mesure de Haar port\'ee par le sous-groupe
compact maximal de $G(\C)$.
\end{lemm}
\begin{proof}
En effet, la derni\`ere formule du lemme~\ref{lemm:norm} implique que le support
de cette mesure est contenu dans le sous groupe compact maximal
de $G(\C)$. D'autre part, les m\'etriques hermitiennes \'etant invariantes
par translation d'un point de torsion, cette mesure est invariante
par translation d'un point de torsion de $G(\C)$. Comme les points de
torsion de $G$ sont denses (pour la topologie complexe)
dans le sous-groupe compact maximal, cette mesure est invariante par
translation par tout point du sous-groupe compact maximal. C'est donc
un multiple de la mesure de Haar.
\end{proof}

Il s'ensuit que l'on peut
remplacer $\Phi$ par $f$ dans le second membre de
l'in\'egalit\'e~\eqref{ineq}.

\begin{lemm} \label{lemm:loinproche}
On a
$\displaystyle \lim_n \int_{O(x_n)} (\Phi-f) \mu_{O(x_n)}=0$.
\end{lemm}
\begin{proof}[Preuve du lemme]
On sait que $h_{\overline{\mathcal L}}(x_n)$ tend vers~$0$.
\'Ecrivons $x_n=(z_n,t_n)$ pour $z_n\in A$ et $t_n\in \gm$, 
de sorte que
\[ h_{\overline{\mathscr L}}(x_n)
    = h_{\overline{\mathscr M}} (z_n)
         + h_{\overline{\mathscr O_\P(1)}} (t_n) \]
et ces deux termes \'etant positifs, ils tendent tous deux vers~$0$.
La d\'ecomposition en facteurs locaux positifs de la hauteur de Weil
sur $\P^1$ implique alors que le terme  local \`a l'infini tend
aussi vers~$0$.

Comme la fonction $(z,t)\mapsto -\log \max(1,\abs t)=\delta(t)$ mesure
la {\og distance\fg} au sous-groupe compact maximal de $G$,
les conjugu\'es de $x_n$ tendent en moyenne
vers ce sous-groupe compact maximal. Pr\'ecis\'ement, si $t\in\gm(\bar\Q)$,
\[ h_{\overline{\mathscr O_\P(1)}} (t)
     \geq - \frac{1}{\# O(t)} \sum_{\tau\in O(t)} \log\max(1,\abs\tau)
     \geq \alpha \frac{\#\{\tau\in O(t)\,;\, \delta(t)\geq\alpha\}}
            {\# O(t)}.
\]
(C'est l'in\'egalit\'e de Bienaym\'e-Tch\'ebitchev.)
Par suite, si $(\alpha_n)$ est une suite de r\'eels positifs,
la proportion des conjugu\'es de $x_n$ dont
la distance au sous-groupe compact maximal de $G(\C)$ est sup\'erieure 
\`a $\alpha_n$ est major\'e par $h(x_n)/\alpha_n$.
D'autre part, la fonction $\Phi-f$ est continue, nulle sur le sous-groupe
compact maximal de $G(\C)$ ; il existe donc une fonction
$\alpha\mapsto\omega(\alpha)$ tendant
vers $0$ en $0$ telle que si la distance de $x$ au sous-groupe compact
maximal est inf\'erieure \`a $\alpha$, $|(f-\Phi)(x)|\leq\omega(\alpha)$.
On a alors l'in\'egalit\'e
\begin{align*}
\abs{ \int_V (\Phi-f) \mu_{O(x_n)} } 
   & \leq  \frac{1}{\# O(x_n)} \sum_{\text{points proches}} 
            + \frac{1}{\# O(x_n)}  \sum_{\text{points loins}} \\
   & \leq  \omega(\alpha_n) + \norm{f-\Phi} \frac{h(x_n)}{\alpha_n} .
\end{align*}
Il suffit de choisir $\alpha_n=\sqrt{h(x_n)}$ pour obtenir le lemme
annonc\'e.

Remarquons que ce lemme n'a pas fait intervenir~$V$ et implique que
toute mesure limite d'une suite de mesures de la forme $(\mu_{O(x_n)})$ 
pour une suite $(x_n)$ de petits points
est concentr\'ee sur le sous-groupe
compact maximal de $G(\C)$.
\end{proof}

En utilisant le lemme, nous obtenons l'in\'egalit\'e
\begin{equation}\tag{$\star\star$}\label{ineq2}
\liminf_n \int_{\bar G} f \mu_{O(x_n)}
\geq \int_G f \nu_{\overline{\mathcal L}}. 
\end{equation}

En appliquant cette in\'egalit\'e \`a $-f$, nous en d\'eduisons le
th\'eor\`eme~\ref{theo:equirep} quand $V=\bar G$,
c'est-\`a-dire le th\'eor\`eme d'\'equir\'epartition
pour une suite g\'en\'erique dans une vari\'et\'e semi-ab\'elienne
de rang torique~$1$.
\end{proof}

\begin{proof}[Preuve du th\'eor\`eme~\ref{theo:equirep} quand $V\neq\bar G$]
Donnons nous une fonction continue $f$ sur $V(\C)$. En vertu des th\'eor\`emes
de Tietze--Urysohn et Stone--Weierstra\ss, nous pouvons supposer
que $f$ est la restriction \`a $V$ d'une fonction $\mathcal C^\infty$
sur $\bar G(\C)$.

Soit $V'\subset V$ l'ensemble des points o\`u la projection
naturelle $\pi:V\ra A$ n'est pas \'etale. Comme $V$
est irr\'eductible, $\dim V' < \dim V$ et son image $\pi(V')$
est de dimension $<\dim A$.

Fixons une famille $(\rho_\lambda)_{\lambda>0}$ 
de fonctions $\mathcal C^\infty$ sur $A(\C)$ telles que pour tout
$\lambda>0$, on a :
\begin{itemize}
\item $0\leq\rho_\lambda \leq 1$ ;
\item $\rho_\lambda=0$ dans un voisinage ouvert $U_\lambda$ de $\pi(V')$ ;
\item $\rho_\lambda=1$ hors d'un voisinage $W_\lambda$ de $\pi(V')$ ;
\end{itemize}
et telle que de plus $\rho_\lambda$ converge simplement
vers la fonction caract\'eristique du compl\'ementaire
dans $\pi(V')$ quand $\lambda\ra 0$.
D\'efinissons ensuite une fonction $f_\lambda$ sur $\bar G(\C)$
en posant $f_\lambda=f \cdot \rho_\lambda\circ\pi$.

\begin{lemm}
Pour tout $\lambda>0$, il existe $\eps_\lambda>0$ tel que pour
$|\eps|<\eps_\lambda$, le courant de type $(1,1)$,
$ c_1(\overline{\mathcal L}|_V) + \eps dd^c f_\lambda $ est positif ou nul.
\end{lemm}
\begin{proof}[Preuve du lemme]
En effet, il est positif dans le voisinage $\pi^{-1}(U_\lambda)\inter V$
de $V'$  puisque $f_\lambda$ y est nul et $c_1(\overline{\mathcal L}|_V)$
y est positif. De plus, sur $V\setminus \pi^{-1}(U_\lambda)$,
le courant 
$ c_1(\overline{\mathcal L}) $
est minor\'e par le courant $\pi^* c_1(\overline{\mathcal M})$,
lequel est strictement positif puisque $\pi:V\ra A$ est \'etale
hors de $\pi^{-1}(U_\lambda)\supset V'$. Ainsi, pour $\eps$
assez petit, $c_1(\overline{\mathcal L}|_V) + \eps dd^c f_\lambda$
est positif hors de $U_\lambda$. 
\end{proof}

Le lemme pr\'ec\'edent montre que nous pouvons appliquer la
proposition~\ref{prop:suz} aux fonctions $f_\lambda$ et $-f_\lambda$, d'o\`u
l'\'egalit\'e
\begin{equation*}\tag{$\star\star\star$}\label{eq3}
\lim_n \int_V f_\lambda \mu_{O(x_n)} = \int_V f_\lambda
\nu_{\overline{\mathcal L}|V}. 
\end{equation*}
Remarquons maintenant que la suite $\pi(x_n)$ est une
suite g\'en\'erique de petits points dans $A$.
D'apr\`es le th\'eor\`eme d'\'equir\'epartition de Szpiro, Ullmo
et Zhang~\cite{szpiro-u-z97},
la suite $\pi(x_n)$ est donc \'equir\'epartie dans $A(\C)$
pour la mesure $\mu_A$ de Haar normalis\'ee de $A(\C)$.
En particulier,
\[  \limsup_n \int_V \abs{f-f_\lambda} \mu_{O(x_n)}
\leq \norm{f} \mu_A( W_\lambda ).  \]
Quand $\lambda\ra 0$, $\mu_A(W_\lambda)$ tend par convergence
domin\'ee vers $\mu_A(V')=0$.

De m\^eme, quand $\lambda\ra 0$, 
$\displaystyle  \int_V f_\lambda \nu_{\overline{\mathcal L}|V}$
converge vers
$\displaystyle \int_V 
    (1-\mathbf{1}_{\pi^{-1}\pi(V')}) f \nu_{\overline{\mathcal L}|V}, $
o\`u $\mathbf 1_{\pi^{-1}\pi(V')}$ d\'esigne la fonction indicatrice de
$\pi^{-1}\pi(V')\subset\bar G$.
Comme $f$ est born\'ee, il suffit donc de montrer que
\[  \int_V \mathbf 1_{\pi^{-1}\pi(V')} 
        c_1(\overline{\mathcal L}_V)^{\dim V} = 0.  \]
Ceci est vrai parce que $\pi^{-1}\pi(V')$
est une sous-vari\'et\'e alg\'ebrique
stricte de $V$, donc un ensemble complet pluripolaire ; un tel
ensemble est ainsi de capacit\'e nulle, c'est-\`a-dire 
n\'egligeable pour des mesures du type Bedford--Taylor
(voir par exemple la proposition~4.6.4 de~\cite{klimek91}, d\'ej\`a
utilis\'ee dans la justification de la d\'efinition~\ref{defi:bt}).


Ainsi, en faisant tendre $\lambda$ vers $0$ dans l'\'egalit\'e~\eqref{eq3},
nous obtenons l'\'egalit\'e du th\'eor\`eme~\ref{theo:equirep},
lequel se trouve donc d\'emontr\'e.
\end{proof}

\section{Preuve de la conjecture de Bogomolov
pour les vari\'et\'es semi-ab\'eliennes isotriviales}
\label{sect:equirep}

\begin{theo}
Si $G$ est une vari\'et\'e ab\'elienne isotriviale,
les conjectures~\ref{conj:bogo}, \ref{conj:pph} et~\ref{conj:equirep}
sont vraies.
\end{theo}

La d\'emonstration de ce th\'eor\`eme est l'objet de ce paragraphe.
Remarquons qu'il nous suffit de d\'emontrer la
conjecture~\ref{conj:equirep} ;
comme il a \'et\'e rappel\'e au paragraphe~2, la conjecture~\ref{conj:bogo}
en est une cons\'equence.
On peut aussi supposer que $G$ est un produit $A\times T$.

La premi\`ere \'etape de la d\'emonstration consiste \`a traiter le cas
d'une vari\'et\'e semi-ab\'elienne dont le rang torique est \'egal \`a~$1$.

\begin{prop}\label{prop:rang1}
Soient $G$ une vari\'et\'e semi-ab\'elienne de rang torique~$1$ sur un
corps de nombres $k$ et
$(x_n)$ une suite de petits points de $G(\bar k)$. On suppose
que cette suite $(x_n)$ est stricte, c'est-\`a-dire qu'aucune sous-suite
n'est contenue dans un sous-groupe alg\'ebrique strict de $G$. Alors, elle
est g\'en\'erique.
\end{prop}
\begin{proof}
Raisonnons par l'absurde et supposons, quitte \`a extraire une sous-suite,
que la suite $(x_n)$ converge vers le point g\'en\'erique d'une sous-vari\'et\'e
irr\'eductible $V$ distincte de $\bar G$.
Comme la hauteur d'un point de $G(\bar\Q)$ est
sup\'erieure ou \'egale \`a celle de son image dans la vari\'et\'e ab\'elienne,
la suite $(\pi(x_n))$ de points de $A(\bar\Q)$ est une suite 
de petits points, stricte puisque la suite $(x_n)$ est stricte dans~$G$.
D'apr\`es le th\'eor\`eme de Zhang~\cite{zhang98}
(conjecture de Bogomolov g\'en\'eralis\'ee
aux sous-vari\'et\'es de vari\'et\'es ab\'eliennes),
$\pi(V)=A$.
Ainsi, ou bien $\dim V=\dim A$  et la projection $V\ra A$ est g\'en\'eriquement
\'etale, ou bien $\dim V>\dim A$ et l'on a $V=\bar G$.

On peut de plus supposer
que le fibr\'e inversible ample $\mathcal L$ choisi
sur $\bar G$ est \'egal \`a $\pi^* \mathcal M\otimes \mathcal O_{\P}(1)$.
Alors, il r\'esulte du th\'eor\`eme~\ref{theo:equirep}
que la suite des mesures de probabilit\'e $(\mu_{O(x_n)})$ 
associ\'ee \`a la suite $(x_n)$ converge
vaguement vers la mesure
\begin{equation}  \label{eq:mesure1}
  \big(\pi^* c_1(\overline{\mathcal M})
            + c_1(\overline{\mathcal O_{\P}(1)}|_V) \big)^{\dim V}
   / \deg_{\mathcal L} V.  
\end{equation}
D'apr\`es le th\'eor\`eme d'\'equir\'epartition de Szpiro, Ullmo et Zhang
(\cite{szpiro-u-z97}), la suite des mesures attach\'ee \`a la suite $\pi(x_n)$ 
converge vaguement pour la mesure
\begin{equation} \label{eq:mesure2}
  \big( c_1(\overline{\mathcal M}) \big)^{\dim A} 
            / \deg_{\mathcal M} A.  
\end{equation}
Cela implique l'\'egalit\'e de l'image directe sur $A$
de la mesure~\eqref{eq:mesure1} et de la mesure~\eqref{eq:mesure2}.
En d\'eveloppant l'\'egalit\'e obtenue,
on obtient, usant du fait que $\pi:V\ra A$ est g\'en\'eriquement fini, 
que les courants 
\[  T_1= c_1(\overline{\mathcal M})^{\dim A-1}\wedge
            \pi_* c_1(\overline{\mathcal O_{\P}(1)}|_V) 
\quad\text{et}\quad
    T_2=c_1(\overline{\mathcal M})^{\dim A}
\]
sont proportionnels.

\begin{lemm}
Le support du courant
$c_1(\overline{\mathcal O_{\P}(1)}|_V)$ est inclus dans l'intersection
de $V(\C)$ et du sous-groupe compact maximal de $G(\C)$.
\end{lemm}
\begin{proof}
Hors de l'intersection de $V(\C)$ et
du sous-groupe compact maximal de $G(\C)$, la formule
donn\'ee dans le lemme~\ref{lemm:norm} montre
que la m\'etrique canonique
sur $\mathcal O_{\P}(1)|_V$ est lisse. On peut donc, sur cet ouvert, calculer
le courant de courbure en restreignant le courant de courbure
de $\mathcal O_{\P}(1)$ sur $\bar G(\C)$. Le calcul du courant a \'et\'e
fait dans le lemme cit\'e et est pr\'ecis\'ement nul hors
du sous-groupe compact maximal.
\end{proof}

Nous reprenons la d\'emonstration de la proposition~\ref{prop:rang1}
en distinguant deux cas.

\medskip
{\itshape\noindent  Premier cas :
supposons que $V\inter D\neq\vide$,}
$D$ d\'esignant le diviseur {\og \`a l'infini\fg}
dans la compactification $\bar G=A\times\P^1$.
Comme dans le lemme~\ref{lemm:norm}, la section canonique $\mathsf s_D$
de $\mathscr O_{\P}(1)$
s'interpr\`ete comme la fonction m\'eromorphe $1/t$.
Comme elle n'est pas constante sur $V$, elle n'est
nulle part localement de module constant sur $V$, si bien qu'on peut trouver
une fibre de $V(\C)\ra A(\C)$ qui ne rencontre pas le sous-groupe compact
maximal de $G(\C)$. En particulier,
$c_1(\overline{\mathcal O_{\P}(1)}|_V)$ est nul au voisinage de cette fibre
et par cons\'equent, le support de $T_1$ est strictement inclus dans $A$.
Comme $T_2$ est une forme diff\'erentielle partout strictement positive sur $A$,
on a n\'ecessairement $T_1=0$.

Or, l'int\'egrale
\[  \int_V \pi^* c_1(\overline{\mathcal M})^{\dim V-1}\wedge
c_1(\overline{\mathcal O_{\P}(1)}|_V)  \]
calcule le degr\'e relativement au fibr\'e inversible ample $\mathcal M$
du cycle $\pi_*(V\cdot D)\in \CH^*(A)$.
La nullit\'e de $T_1$
implique que ce degr\'e est nul.
Or, l'intersection $V\cdot D$ est propre,
$V$ n'\'etant pas contenue dans $D$, si bien qu'un repr\'esentant
de $V\cdot D$ est donn\'e par une somme $\sum_i n_i [V_i]$, o\`u
les $V_i$ sont les composantes irr\'eductibles de $V\inter D$
et les $n_i$ des entiers strictement positifs.
On obtient que $V\inter D=\vide$, ce qui est absurde.

\medskip
{\itshape\noindent Deuxi\`eme cas :
supposons que $V\inter D=\vide$.}
Dans ce cas, $V$ est affine et propre sur $A$ et $\pi:V\ra A$
est fini.
Cela fournit
une \'equation pour $V$ dans $A\times\P^1$ :
\[  t^d + a_1(z) t^{d-1}+ \cdots +a_d(z) =0.  \]
Comme $V\inter D=\vide$, on voit que pour tout $i$,
$a_i$ est une fonction r\'eguli\`ere sur $A$, donc une constante.
Autrement dit, $V=A\times V_1$, o\`u $V_1$ est un ensemble \textit{fini}
de points de ${\mathbf A}^1$. On en d\'eduit qu'il existe une suite
de points de $V_1$ dont les hauteurs tendent vers~$0$, si bien
que $V_1$ contient un point de torsion. Puisque $V$ est irr\'eductible,
$V_1$ est form\'e de points de torsion et $V$ est une sous-vari\'et\'e de torsion,
contrairement \`a l'hypoth\`ese que la suite $(x_n)$ est stricte.
Nous avons donc d\'emontr\'e que la suite $(x_n)$ est g\'en\'erique dans $G$.
\end{proof}


\begin{proof}[Preuve de la conjecture~\ref{conj:equirep}]
Soit $(x_n)$ une suite stricte de petits points de $G(\bar k)$.
Pour simplifier les notations, on note $\nu_{G}$ la mesure
de Haar normalis\'ee port\'e par le sous-groupe compact maximal
de $G(\C)$.

Pour tout caract\`ere du tore $\psi:T\ra\gm$, on peut pousser
la situation le long de $\psi$, ce qui nous fournit
une vari\'et\'e semi-ab\'elienne $G_\psi=A\times \gm$ de rang torique~$1$,
munie d'un morphisme $\psi:G\ra G_\psi$, 
ainsi qu'une suite de petits points $x_{\psi,n}=\psi(x_n)\in G_\psi(\bar k)$.
Sauf si $\psi$ est le caract\`ere trivial, cette suite est stricte.
D'apr\`es la proposition pr\'ec\'edente, elle est donc g\'en\'erique.

La suite des orbites $\delta_{O(x_{\psi,n})}$ est, en vertu
de la proposition~\ref{prop:rang1}, \'equir\'epartie pour
la mesure de Haar normalis\'ee $\nu_{G_\psi}$.

Alors, si $\mu$ est une valeur d'adh\'erence de la suite de mesures 
$\delta_{O(x_n)}$, on voit que pour tout caract\`ere $\psi$ non trivial,
$\psi_*\mu = \nu_{G_\psi}$.

D'autre part, le support de $\mu$ est inclus dans l'intersection
des images r\'eciproques par $\psi$ du sous-groupe compact maximal
de $G_\psi(\C)$, c'est-\`a-dire que $\mu$ est une mesure 
de probabilit\'e port\'ee
par le sous-groupe compact maximal de $G(\C)$.

Or, la mesure $\nu_G$ a toutes ces propri\'et\'e 
et la th\'eorie des s\'eries de Fourier implique
que $\mu=\nu_G$.
Comme l'ensemble des mesures de probabilit\'e sur $\bar G(\C)$ est
compact pour la topologie faible-$*$ (th\'eor\`eme de Banach--Alaoglu),
la suite $\delta_{O(x_n)}$ converge vers la mesure $\nu_G$,
ainsi qu'il fallait d\'emontrer.
\end{proof}

Nous avons donc prouv\'e l'analogue
de la conjecture de Bogomolov pour une vari\'et\'e
semi-ab\'elienne isotriviale.

\section*{Appendice : Le th\'eor\`eme d'\'equir\'epartition de Bilu via la g\'eom\'etrie d'Arakelov}

Nous d\'emontrons dans ce paragraphe un th\'eor\`eme d'\'equir\'epartition
pour une suite g\'en\'erique de petits points
(ce qui revient \`a stricte, dans ce cas) dans $\P^1$ d\^u initialement \`a
Bilu~\cite{bilu97}. Il a montr\'e comment en d\'eduire la conjecture
de Bogomolov pour un tore ({\em i.e.}~quand $A=0$, avec nos notations).

\begin{theo*}[Bilu] \label{bilu}
Soit $(x_n)_{n\in\N}$ une suite d'\'el\'ements de $\gm(\bar\Q)=\bar\Q^*$.
On suppose~:
\begin{itemize}
\item qu'aucune sous-suite n'est constante ({\em i.e.,} la suite $(x_n)$ est
g\'en\'erique)~;
\item que la hauteur normalis\'ee de $x_n$ tend vers~$0$ quand
$n\ra+\infty$.
\end{itemize}
Notons $O(x_n)$ l'adh\'erence
de $x_n\in\P^1(\bar\Q)$ dans $\P^1_\Q$.
Alors, pour toute fonction continue $f:\P^1(\C)\ra\R$, on a la
convergence~:
\[  \frac{1}{\Card O(x_n)} \sum_{x\in O(x_n)(\C)} f(x) \rightarrow 
   \frac{1}{2\pi} \int_{0}^{2\pi} f(e^{i\theta})\, d\theta.  \]
Autrement dit la suite des mesures $\frac{1}{\Card O(x_n)}
\delta_{O(x_n)(\C)}$ converge vaguement vers la mesure de Haar normalis\'ee
sur le sous-groupe compact maximal de $\gm$.
\end{theo*}
\begin{proof}
La densit\'e des fonctions de classe $\mathscr C^ \infty$ dans les
fonctions continues sur $\P^1(\C)$ nous permet de supposer
que $f$ est $\mathcal C^\infty$. Soit alors $\Phi:\P^1(\C)\ra\R$
continue telle que~:
\begin{itemize}
\item $\Phi=f$ sur $S^1$~;
\item $\Phi$ est harmonique sur les deux disques de $\P^1(\C)$ bord\'es
par $S^1$ (d'\'equations en coordonn\'ee inhomog\`ene, $|z|<1$ et $|z|>1$).
\end{itemize}
L'existence et l'unicit\'e de $\Phi$ est classique (noyau de Poisson),
on v\'erifie par un calcul explicite que $dd^c\Phi$ est une mesure
port\'ee par le cercle unit\'e,
absolument continue par rapport \`a la mesure de Haar sur ce cercle.
En effet, soit $f(e^{i\theta})=\sum_{n\in\Z}c_n e^{in\theta}$ la s\'erie
de Fourier de $f$ sur le cercle unit\'e ; on pose
\[  \Phi(r e^{i\theta}) = \sum_{n\in\Z} c_n \min(r^n,r^{-n}) e^{in\theta}.
\] 
Alors,
\[  dd^c\Phi = - \left(\sum_{n\in\Z} |n|c_n e^{in\theta}\right) \,
\frac{d\theta}{2\pi}\wedge \delta_{r=1},  \]
est une mesure concentr\'ee sur le cercle unit\'e, absolument continue
par rapport \`a la mesure de Haar sur le cercle unit\'e.

Soit $\overline{\mathcal L}$ le fibr\'e inversible hermitien $\mathcal
O_{\P^1}(1)$ muni de sa m\'etrique canonique, dont la forme
de courbure est la mesure de Haar $\displaystyle \frac{1}{2\pi} {d\theta}$ 
port\'ee par le sous-groupe compact maximal de $\C^*$.
Pour $\lambda\in\R$,
d\'efinissons $\overline{\mathcal L}(\lambda\Phi)$ comme le fibr\'e en
droites hermitien suivant : c'est $\mathcal L$ comme fibr\'e en droites,
la m\'etrique \'etant celle de $\overline{\mathcal L}$ multipli\'ee
par $\exp(-\lambda\Phi)$.
La forme de courbure  de $\overline{\mathcal L}(\lambda\Phi)$ est
\'egale \`a
\[  c_1(\overline{\mathcal L}) + \lambda dd^c\Phi .  \]
Autrement dit, elle est nulle hors du sous-groupe compact maximal
et positive sur le cercle unit\'e si $|\lambda|$ est assez petit.

D'apr\`es le th\'eor\`eme~1.10 de~\cite{zhang95b}
sur l'amplitude arithm\'etique,
cons\'equence du th\'eor\`eme
de Hilbert--Samuel arit\-hm\'e\-tique,
on a
\[  \liminf_n h_{\overline{\mathcal L}(\lambda\Phi)}(x_n) \geq
h_{\overline{\mathcal L}(\lambda\Phi)}(\P^1),  \]
soit
\[  \liminf_n \lambda \left( \frac{1}{\Card O(x_n)} \sum_{x\in O(x_n)}
\Phi(x)\right)
    \geq  \lambda \int_{\P^1(\C)}\Phi c_1(\overline{\mathcal L})
+O(\lambda^2).  \]
Si l'on fait tendre $\lambda $ vers $0$, par valeurs sup\'erieures ou
inf\'erieures, on obtient alors
\[  \lim_n \frac{1}{\Card O(x_n)} \sum_{x\in O(x_n)}
\Phi(x) = \frac{1}{2\pi} \int_{0}^{2\pi} \Phi(e^{i\theta})\, d\theta.  \]

Comme $\Phi=f$ sur $S^1$, on peut remplacer $\Phi$ par $f$ dans
cette derni\`ere int\'egrale.
D'autre part, la hauteur d'un point $\xi\in\bar\Q^*$ est minor\'ee
par la demi-somme des composantes archim\'ediennes des hauteurs
de $\xi$ et $\xi^{-1}$, c'est-\`a-dire
\[  h(\xi) \geq  \frac{1}{2[\Q(\xi):\Q]} 
    \sum_{\sigma:\Q(\xi)\ra\C} \log \max ( |\sigma(\xi)|, |\sigma(\xi)|^{-1}),  \]
soit
\[  h(\xi) = \frac{1}{\Card O(\xi)} \sum_{x\in O(\xi)} \frac12\log \max
(|x|,|x|^{-1}).  \]
Notons que la fonction 
$x\mapsto d(x)=\frac12\log\max(|x|,|x|^{-1})$ mesure la {\og distance\fg} 
de $x$ au cercle unit\'e. 

Comme la suite $h(x_n)$ tend vers~$0$ et comme aucune sous-suite n'est
constante,
le th\'eor\`eme de Northcott implique que 
$\lim \Card O(x_n)=+\infty$. Soit pour tout $n$
un r\'eel $\alpha_n>0$ ; on suppose que $\alpha_n\ra 0$. 
La proportion des $x\in O(x_n)$ tels
que $d(x)\geq\alpha_n$ s'estime ainsi~:
\[  \Card \{ x\in O(x_n)\sozat d(x)\geq\alpha_n\}
    \leq \Card O(x_n) \frac{h(x_n)}{\alpha_n}.  \]
Remarquons alors que
\begin{align*}
\abs{\sum_{x\in O(x_n)} f(x)-\sum_{x\in O(x_n)}\Phi(x) } 
&\leq 
   \sum_{d(x)\geq\alpha_n} \abs{f(x)-\Phi(x)} + \sum_{d(x)\leq\alpha_n}
\abs{f(x)-\Phi(x)} \\
& \leq 
   \Card\{x\in O(x_n)\sozat d(x)\geq\alpha_n\}  \norm{f-\Phi}_\infty \\
& \qquad
  {} + \Card O(x_n) \sup_{d(x)\leq\alpha_n} |f-\Phi| \\
& \leq 
   \Card O(x_n) \left(  \frac{h(x_n)}{\alpha_n} \norm{f-\Phi}_\infty
  + o(1) \right)
\end{align*}
d'apr\`es l'uniforme continuit\'e de $f-\Phi$ sur $\P^1(\C)$
et le fait que $\alpha_n\ra 0$. On choisit maintenant
$\alpha_n=\sqrt{h(x_n)}$. Ainsi, 
\[  \frac{1}{\Card O(x_n)} \sum_{x\in O(x_n)} f(x)
- \frac{1}{\Card O(x_n)} \sum_{x\in O(x_n)}\Phi(x)  \]
converge vers~$0$.

En d\'efinitive, nous avons prouv\'e que
\[  \frac{1}{\Card O(x_n)} \sum_{x\in O(x_n)} f(x)  \]
converge vers
\[  \frac{1}{2\pi} \int_{0}^{2\pi} f(e^{i\theta})\, d\theta,  \]
ce qui conclut la preuve du th\'eor\`eme.
\end{proof}

\def\noop#1{\ignorespaces}
\bibliographystyle{smfplain}
\bibliography{acl}

\end{document}